\font\tta=cmr10 scaled\magstep1

\font\ttc=cmss12

\font\tte=cmbsy10
\font\ttcx=cmss8

\font\ttf=cmr8

\font\tth=cmbx10 scaled\magstep2

\def\ov#1{\overline{#1}}
\def\projx#1{\hbox{{\tta \hskip 4pt P\hskip -10pt I\hskip 6pt}$(#1)$}}
\def\proj#1{\hbox{{\tta \hskip 2pt P\hskip -10pt I~$_#1$}}}

\def\cmplxb{\ \hbox{{\ttc \hskip -2pt C\hskip -5pt I \hskip -1pt}}}
\def\cx#1{\cmplxb^{#1}}
\def\scx#1{ \hbox{\ttcx  C\hskip -4pt l }^{#1}}
\def\codim{\mathop{\rm codim}\nolimits}
\def\coker{\mathop{\rm coker}\nolimits}
\def\Hom{\mathop{\rm Hom}\nolimits}
\def\Ext{\mathop{\rm Ext}\nolimits}
\def\imm{\mathop{\rm Im}\nolimits}

\def\lra{\longrightarrow}
\def\som{\mathop{\hbox{$\displaystyle\bigoplus$}}\limits}
\def\sigg{\mathop{\hbox{$\displaystyle\sum$}}\limits}

\def\supp{\mathop{\hbox{$\sup$}}\limits}
\def\carre{$\Box$}
\def\timex{\setbox250=\hbox{$\times$}\hskip 5pt\Box\hskip -0.75em
{\raise 1.5pt\vbox{\box250}}\hskip 5pt}
\def\paragra{{\tte \char120}}
\def\para{\paragra~\hskip -2pt}
\def\subsubsubsection#1#2{\noindent{\bf #1} \ \ {\it #2}}
\def\hfl#1#2{\smash{\mathop{\hbox to 12mm{\rightarrowfill}}
\limits^{\scriptstyle#1}_{\scriptstyle#2}}}
\def\vfl#1#2{\llap{$\scriptstyle #1$}\left\downarrow
\vbox to 6mm{}\right.\rlap{$\scriptstyle #2$}}

\def\diagram#1{\def\normalbaselines{\baselineskip=0pt
\lineskip=10pt\lineskiplimit=1pt}  \matrix{#1}}
\def\pline#1{<\hskip-3.5pt#1\hskip-3.5pt>}
\def\q#1#2{{{#1}\over{#2}}}
\def\m#1{{\hbox{$#1$}}}
\def\fexc{fibr\'e exceptionnel }
\def\fexcs{fibr\'es exceptionnels }
\def\ot{\otimes}
\def\E{{\cal E}}
\def\F{{\cal F}}
\def\G{{\cal G}}
\def\O{{\cal O}}
\def\og{\leavevmode\raise.3ex\hbox{$\scriptscriptstyle\langle\!\langle$}}
\def\fg{\leavevmode\raise.3ex\hbox{$\scriptscriptstyle\,\rangle\!\rangle$}}
\def\dem{\noindent{\em D\'emonstration}. }

\def\Oad{On a un diagramme commutatif avec lignes et colonnes exactes :}
\def\g#1{\ifnum #1=1 \gamma_1 \else {
\ifnum #1=2 \gamma_2 \else {
\ifnum #1=3 \gamma_3 \else {
\ifnum #1=4 \gamma_4 \else {}
\fi}\fi}\fi}\fi}
\def\X#1{\ifnum #1=1 X_1 \else {\ifnum #1=2 X_2 \else{
\ifnum #1=3 X_3 \else{ \ifnum #1=4 X_4 \else{}
\fi}\fi}\fi}\fi}
\def\x#1{\ifnum #1=3 x_3 \else{\ifnum #1=4 x_4 \else{} \fi}\fi}
\def\GG#1{\ifnum #1=0 G_0 \else{\ifnum #1=1 G_1 \else{\ifnum #1=2 G_2
\else {}
\fi}\fi}\fi}
\def\gg#1{\ifnum #1=0 g_0 \else{\ifnum #1=1 g_1 \else{\ifnum #1=2 g_2
\else {}
\fi}\fi}\fi}
\def\psix{\phi_2}
\documentstyle[twoside,leqno,12pt]{article}
\begin{document}
\newtheorem{xprop}{Proposition}[section]
\newtheorem{xlemm}[xprop]{Lemme}
\newtheorem{xtheo}[xprop]{Th\'eor\`eme}
\newtheorem{xcoro}[xprop]{Corollaire}
\newtheorem{quest}{Question}
\newtheorem{defin}{D\'efinition}

\markboth{Jean-Marc Dr\'ezet}{Espaces abstraits de morphismes et
mutations}

\title{Espaces abstraits de morphismes et mutations}
\author{Jean-Marc Dr\'ezet}
\date{}
\maketitle

\def\refname{R\'ef\'erences}
{\ttf
\noindent  \ Universit\'e Paris 7, UMR 9994 du CNRS,
Aile 45-55, 5\m{^e} \'etage

\vskip -0.2cm

\noindent \ 2, place Jussieu, F-75251 Paris Cedex 05, France

\vskip -0.2cm

\noindent \ e-mail : drezet@mathp7.jussieu.fr
}

\bigskip

\bigskip

\bigskip

{\parindent 3cm
{\tth \hskip 3cm Sommaire}

\bigskip

1 - Introduction

\medskip

2 - Un exemple simple

\medskip

3 - Vari\'et\'es de modules de morphismes de type $(r,s)$

\medskip

4 - Mutations d\'efinies \`a l'aide de la suite spectrale de Beilinson

\hskip 0.7cm g\'en\'eralis\'ee

\medskip

5 - Mutations en termes de morphismes de faisceaux

\medskip

6 - Espaces abstraits de morphismes

\medskip

7 - Mutations de morphismes de type \m{(r,s)}

\medskip

8 - Applications

}

\bigskip

\bigskip

\section{Introduction}

\subsection{Vari\'et\'es de modules de morphismes}

Soient \m{X} une vari\'et\'e alg\'egrique projective sur le corps
des nombres complexes, et \m{\E}, \m{\F}  des faisceaux
alg\'ebriques coh\'erents sur \m{X}. Soit
$$W = \Hom(\E,\F) .$$
Alors le groupe alg\'ebrique
$$G = Aut(\E)\times Aut(\F)$$
agit d'une fa\c con \'evidente sur \m{W}. Si deux morphismes sont dans
la m\^eme \m{G}-orbite, leurs noyaux sont isomorphes, ainsi que
leurs conoyaux. C'est pourquoi il peut \^etre int\'eressant, pour
d\'ecrire certaines vari\'et\'es de modules de faisceaux, de
construire de bons quotients d'ouverts \m{G}-invariants de \m{W} par
\m{G}. On s'int\'eresse au cas particulier suivant : soient \m{r}, \m{s} des
entiers positifs, \m{\E_1,\ldots,\E_r,},
\m{\F_1,\ldots,\F_s} des faisceaux
coh\'erents sur \m{X}, qui sont {\em simples}, c'est-\`a-dire que leurs
seuls endomorphismes sont les homoth\'eties. On suppose aussi que
$$\Hom(\E_i,\E_{i'}) = \lbrace 0 \rbrace \ \ {\rm si \ } i > i' \ , \
\Hom(\F_j,\F_{j'}) = \lbrace 0 \rbrace \ \ {\rm si \ } j > j', $$
$$\Hom(\F_j,\E_i) = \lbrace 0 \rbrace \ \
{\rm pour \ tous \ } i,j .$$
Soient \m{M_1,\ldots,M_r}, \m{N_1,\ldots,N_s} des espaces vectoriels
complexes de dimension finie. On suppose que
$$\E = \som_{1\leq i\leq r}(\E_i\ot M_i) \ , \
\F = \som_{1\leq l\leq s}(\F_l\ot N_l) .$$
Les \'el\'ements de \m{W} sont appel\'es {\em morphismes de type} \m{(r,s)}.
Le groupe \m{G} n'est pas r\'eductif en g\'en\'eral.
On a consid\'er\'e dans
\cite{dr_tr} le probl\`eme de l'existence de bon quotients d'ouverts
\m{G}-invariants de \m{\Hom(\E,\F)}.
On introduit une notion de {\em semi-stabilit\'e} pour les
mor-\break phismes de type \m{(r,s)} qui d\'epend du choix d'une suite
\m{(\lambda_1,\ldots,\lambda_r,}\m{\mu_1,\ldots,\mu_s)} de nombres
rationnels positifs tels que
$$\sigg_{1\leq i\leq r}\lambda_i\dim(M_i) =
\sigg_{1\leq l\leq s}\mu_l\dim(N_l) = 1.$$
On appelle cette suite une {\em polarisation} de l'action de \m{G}.
Il existe un bon quotient de l'ouvert des points
semi-stables pour certaines valeurs de
\m{(\lambda_1,\ldots,\lambda_r,} \m{\mu_1,\ldots,\mu_s)}
(ces r\'esultats sont rappel\'es au \para 3).

Pour traiter ce genre de probl\`eme de la mani\`ere la plus g\'en\'erale,
on le traduit d'abord en termes
d'alg\`ebre lin\'eaire (c'est-\`a-dire que ce qu'on \'etudie est une action
particuli\`ere d'un certain groupe non r\'eductif sur un espace vectoriel de
dimension finie). C'est ce qui est fait dans \cite{dr_tr} et ici.

Les morphismes de type \m{(2,1)} sont utilis\'es dans \cite{dr2} pour
d\'ecrire certaines vari\'et\'es de modules de faisceaux semi-stables
sur \proj{2}. Dans un certain nombre de travaux (cf. par exemple
\cite{miro}, \cite{oko}) des faisceaux semi-stables ou des faisceaux
d'id\'eaux de sous-vari\'et\'es de l'espace projectif sont d\'ecrits
comme conoyaux de morphismes de type \m{(r,s)}.

\subsection{Mutations de morphismes}

Le but du pr\'esent article est de d\'ecrire et d'\'etudier certaines
transformations, appel\'ees {\em mutations}, associant \`a un morphisme de
type \m{(r,s)} un autre morphisme, pouvant \^etre d'un autre type
(mais la somme \m{r+s} reste constante). On obtient en quelque sorte
une correspondance entre deux espaces de morphismes \m{W} et \m{W'}, sur
lesquels agissent respectivement les groupes en g\'en\'eral non
r\'eductifs \m{G} et \m{G'}. Ceci permet
de d\'efinir une bijection de l'ensemble des \m{G}-orbites
d'un ouvert de \m{W} sur l'ensemble des \m{G'}-orbites d'un ouvert de
\m{W'}. La forme que prend une mutation dans le language des morphismes
de faisceaux est explicit\'ee au \para 1.3 (dans la description du
chapitre 5). On donnera toutefois une d\'efinition plus abstraite et
plus g\'en\'erale de ce qu'est une mutation dans le chapitre 6.

On associe de mani\`ere naturelle
\`a chaque polarisation \m{\sigma} de l'action de \m{G} sur \m{W} une
polarisation \m{\sigma'} de l'action de \m{G'} sur \m{W'}.
Dans certains cas on montre qu'un point de \m{W} est semi-stable
relativement \`a \m{\sigma} si et seulement si la mutation de ce
point est semi-stable relativement \`a \m{\sigma'}. Ceci
permet de prouver que les quotients correspondants sont isomorphes. On
peut ainsi \'etendre les r\'esultats de \cite{dr_tr} \`a d'autres
polarisations. Par exemple, dans \cite{dr2} on prouve dans certains cas
l'existence de bon quotients (lesquels sont isomorphes \`a des vari\'et\'es
de modules de faisceaux semi-stables sur \m{\proj{2}}). On
n'a pas besoin dans ce cas d'un th\'eor\`eme d'existence d'un quotient par
un groupe non r\'eductif, car la vari\'et\'e de modules existe d\'ej\`a.
Ces exemples de bons quotients ne peuvent pas \^etre directement retrouv\'es
\`a partir de \cite{dr_tr}, mais en utilisant des mutations, on peut se
ramener aux cas trait\'es dans \cite{dr_tr}. Autre exemple, les morphismes
$$\O(-2)\oplus\O(-1)\lra\O\ot\cx{n+2}$$
sur \m{\proj{n}}. L'application directe de \cite{dr_tr} ne fournit
aucun quotient non vide. En utilisant des mutations,
on peut trouver plusieurs types de quotients.

La d\'efinition des mutations s'introduit naturellement lorsqu'on \'etudie
les faisceaux semi-stables sur \m{\proj{n}} au moyen des
suites spectrales de Beilinson g\'en\'eralis\'ees. On associe une telle suite
spectrale \`a un faisceau coh\'erent \m{\E}
sur \m{\proj{n}} et \`a une {\em base d'h\'elice \m{\sigma}
de fibr\'es exceptionnels } sur \m{\proj{n}} (cf \cite{dr1} et
\cite{go_ru}). Si le diagramme de Beilinson correspondant est suffisamment
simple, on obtient une suite exacte
$$0\lra \som_{1\leq i\leq r}(E_i\otimes\cx{m_i})
\lra\som_{1\leq l\leq s}(F_l\otimes\cx{n_l})\lra\E\lra 0 ,$$
la base d'h\'elice \m{\sigma} \'etant
\m{(E_1,\ldots,E_r,F_1,\ldots,F_s)} (et donc \m{r + s = n+1}). On peut,
en changeant judicieusement la base d'h\'elice, obtenir d'autres
repr\'esentations semblables de \m{\E}. On peut changer de
base d'h\'elice en faisant subir \`a celle dont on part une s\'erie de
transformations \'el\'ementaires appel\'ees {\em mutations}, d'o\`u la
terminologie employ\'ee pour les transformations de morphismes \'etudi\'ees
ici.
Pour la d\'efinition, les propri\'et\'es et
l'usage des suites spectrales de Beilinson
g\'en\'eralis\'ees sur les espaces
projectifs, voir \cite{dr1}, \cite{dr2}, \cite{dr3}, \cite{dr_lp},
\cite{go_ru}.

\subsection{Plan des chapitres suivants}

Dans le chapitre 2 on donne un exemple simple et bien connu de mutations
dans le cas des morphismes de type (1,1). C'est ce type de r\'esultats
qu'il s'agit de g\'en\'eraliser.

\medskip

Dans le chapitre 3 on rappelle certains r\'esultats de \cite{dr_tr},
concernant les quotients d'espa-\break ces de morphismes de type \m{(r,s)}.
Le th\'eor\`eme 3.1 d\'ecrit ce qu'on sait des quotients d'espaces
de morphismes de type (2,1).

\medskip

Dans le chapitre 4 on rappelle la d\'efinition des suites spectrales de
Beilinson g\'en\'erali-\break s\'ees sur les espaces projectifs et on
d\'ecrit des mutations de morphismes de type
\m{(r,s)} obtenues en utilisant les suites spectrales de
Beilinson g\'en\'eralis\'ees sur les espaces projectifs.

\medskip

Dans le chapitre 5 on g\'en\'eralise un peu ce qui pr\'ec\`ede.
On donne la d\'efinition des mutations de morphismes
en termes de faisceaux. Plus pr\'ecis\'ement on montre que si un
faisceau coh\'erent peut \^etre repr\'esent\'e comme conoyau d'un
morphisme injectif de faisceaux, on peut dans certaines conditions le
repr\'esenter aussi comme conoyau d'un morphisme injectif d'un autre
type. Les r\'esultats du \para 5.1 sont plus g\'en\'eraux que ce qui est
n\'ecessaire ici. Dans le \para 5.2 on donne des applications aux morphismes
de type \m{(r,s)}. Voici un exemple du type de r\'esultat obtenu : soit
$$\Phi  : \som_{1\leq i\leq r}(\E_i\ot M_i) \lra
\som_{1\leq l\leq s}(\F_l\ot N_l)$$
un morphisme injectif, \m{\cal U} son conoyau
et \m{p} un entier tel que \m{0\leq p\leq r-1}. On suppose que pour
\m{p+1\leq j\leq r} le morphisme canonique
$$\E_j\lra \Hom(\E_j,\F_1)^*\ot \F_1$$
est injectif. Soit \m{\G_j} son conoyau. Soit
$$f_p : \som_{p+1\leq j\leq r}(\Hom(\E_j,\F_1)^*\ot M_j)\lra N_1$$
l'application lin\'eaire d\'eduite de $\Phi$. On suppose que $f_p$ est
surjective. Alors on montre que sous certaines hypoth\`eses
il existe une suite exacte
$$0\lra \biggl(\bigoplus_{1\leq i\leq p}(\E_i\ot M_i)\biggr)\oplus
(\F_1\ot\ker(f_p))\lra
\biggr(\bigoplus_{p<j\leq r}(\G_j\ot M_j)\biggl)\oplus
\biggr(\bigoplus_{2\leq l\leq s}(\F_l\ot N_l)\biggl)\lra{\cal U}\lra 0.$$
On a donc associ\'e \`a un morphisme de type \m{(r,s)} un morphisme de
type

\noindent\m{(p+1,r+s-p-1)}.

\medskip

Dans le chapitre 6 on d\'efinit les mutations dans un cadre plus abstrait.
On d\'efinit des actions de groupes sur des espaces vectoriels
model\'ees sur les cas \'etudi\'es dans le chapitre pr\'ec\'edent.
Une mutation est dans ce cas une correspondance entre une telle action d'un
groupe $G$ sur un espace vectoriel $V$ et une autre action d'un groupe $G'$
sur un espace vectoriel $V'$, de telle sorte qu'on ait une bijection
\ \m{V^0/G\simeq {V'}^0/G'}, pour des ouverts ad\'equats $V^0$ et
\m{{V'}^0} non vides de $V$ et $V'$ respectivement.

\medskip

Dans le chapitre 7, on applique les r\'esultats qui pr\'ec\`edent dans le
but de trouver d'autres cas o\`u on sait d\'efinir des bons quotients
d'espaces de morphismes de type \m{(r,s)}. Dans le cas des morphismes
de type (2,1), le th\'eor\`eme 7.6 obtenu \'etend les r\'esultats du
th\'eor\`eme 3.1.

\medskip

Dans le chapitre 8 on donne des exemples d'applications des r\'esultats
pr\'ec\'edents.

\bigskip

{\bf Remerciements.} L'auteur tient \`a remercier G. Trautmann pour de
nombreuses discussions
qui l'ont beaucoup aid\'e, ainsi que l'Universit\'e de Kaiserslautern pour
son hospitalit\'e durant la r\'ealisation d'une partie de ce travail.

\vfil
\eject

\section{Un exemple simple}

Les r\'esultats de ce chapitre sont d\'emontr\'es dans \cite{dr2}. Soient
\m{L}, \m{M} et \m{N} des espaces vectoriels complexes de dimension finie,
avec \ \m{\dim(L)\geq 3}. On pose \
\m{q = \dim(L)}, \m{m = \dim(M)}, \m{n = \dim(N)}.
Les applications lin\'eaires
$$L\ot M\lra N$$
sont appel\'ees des \m{L}-{\em modules de Kronecker}. Soit
$$W=\Hom(L\ot M,N).$$
Sur \m{W} op\`ere de mani\`ere \'evidente le groupe alg\'ebrique r\'eductif
$$G=(GL(M)\times GL(N))/\cx{*}.$$
L'action de \ \m{SL(M)\times SL(N)} \ sur \m{\projx{W}} se lin\'earisant de
fa\c con \'evidente, on a une notion de point {\em (semi-)stable} de
\m{\projx{W}} (au sens de la g\'eom\'etrie invariante). On montre que si
\m{f\in W}, \m{f} est semi-stable (resp.
stable) si et seulement si pour tous sous-espaces vectoriels \m{M'} de
\m{M} et \m{N'} de \m{N}, tels que \m{M'\not = \lbrace 0\rbrace},
\m{N'\not =N}, et \ \m{f(L\ot M')\subset N'}, on a
$$\q{\dim(N')}{\dim(M')}\geq\q{\dim(N)}{\dim(M)}\ \ {\rm (resp.}
\ > \ {\rm)}.$$
Soit \m{W^{ss}} (resp. \m{W^s}) l'ouvert des points semi-stables (resp.
stables) de \m{W}. Alors il existe un bon quotient (resp. un quotient
g\'eom\'etrique)
$$N(L,M,N) = W^{ss}//G \ \ \ {\rm (resp. } \ \ N_s(L,M,N)=W^s/G \ {\rm)},$$
\m{N(L,M,N)} est projective, et \m{N_s(L,M,N)} est un ouvert lisse de
\m{N(L,M,N)}.

On pose \ \m{m'=qm-n}. On suppose que \m{m'>0}. Soit \m{M'} un
espace vectoriel complexe de dimension \m{m'}. Soit \ \m{f:L\ot M\lra N} \
un \m{L}-module de Kronecker surjectif. Alors \ \m{\dim(\ker(f))=m'}. Soient
$$f' : L^*\ot\ker(f)\lra H_0$$
la restriction de l'application
$$tr\ot I_{H_0} : L^*\ot L\ot H_0\lra H_0$$
(\m{tr} d\'esignant l'application trace), et
$$A(f) : L^*\ot H_0^*\lra\ker(f)^*$$
l'application lin\'eaire d\'eduite de \m{f'}, qu'on peut voir comme un
\'el\'ement de

\noindent\m{W' = \Hom(L^*\ot H_0^*,M')}, en utilisant un
isomorphisme \ \m{\ker(f)^*\simeq M'}. Soit \m{W_0} l'ouvert de \m{W}
constitu\'e des applications surjectives, \m{W'_0} l'ouvert analogue de
\m{W'}, et

\noindent\m{G' = (GL(H_0^*)\times GL(M'))/\cx{*}}. On d\'emontre
ais\'ement la

\begin{xprop}
1 - En associant \m{A(f)} \`a \m{f} on d\'efinit une bijection
$$W_0/G \ \simeq \ W'_0/G'.$$

\noindent 2 - La bijection pr\'ec\'edente induit un isomorphisme
$$N(L,M,N)\ \simeq \ N(L^*,H_0^*,M')$$
(induisant un isomorphisme \ \m{N_s(L,M,N)\ \simeq \ N_s(L^*,H_0^*,M')}).
\end{xprop}

\section{Vari\'et\'es de modules de morphismes de type $(r,s)$}

On rappelle ici le probl\`eme des vari\'et\'es de modules de morphismes
de type $(r,s)$ abord\'e dans \cite{dr_tr}.

En termes de faisceaux, on consid\`ere des morphismes
$$\E = \som_{1\leq i\leq r}(\E_i\ot M_i) \lra
\som_{1\leq l\leq s}(\F_l\ot N_l) = \F ,$$
et l'action du groupe alg\'ebrique
$$G = Aut(\E)\times Aut(\F)$$
\ sur l'espace vectoriel \m{W} de tous ces morphismes.

Il est pr\'ef\'erable de g\'en\'eraliser ce probl\`eme en termes
d'alg\`ebre lin\'eaire.

\subsection{D\'efinition abstraite de \m{W}}

Soient \m{r,s} des entiers positifs,
\m{H_{li}}, \m{A_{ji}}, \m{B_{ml}}, \m{1\leq i\leq j\leq r},
\m{1\leq l\leq m\leq s} des espaces vectoriels de dimension finie (qui
jouent le r\^ole de \m{\Hom(\E_i,\F_l)},
\m{\Hom(\E_i,\E_j)} et \m{\Hom(\F_l,\F_m)}
respectivement). On suppose que  \ \m{A_{ii} = \cx{}}
\ pour \m{1\leq i\leq r} et \ \m{B_{ll} = \cx{}} \
pour \m{1\leq l\leq s}. Pour \m{1\leq i\leq j\leq k\leq r} et
\m{1\leq l\leq m\leq n\leq s} on se donne des applications lin\'eaires
(appel\'ees {\it compositions})
$$H_{lj}\ot A_{ji}\lra H_{li},$$
$$A_{kj}\ot A_{ji}\lra A_{ki},$$
$$B_{ml}\ot H_{li}\lra H_{mi},$$
$$B_{nm}\ot B_{ml}\lra B_{nl}.$$
On suppose que si \m{i=j} les deux premi\`eres applications sont les
identit\'es, ainsi que la se-\break conde si \m{j=k}, la quatri\`eme si
\m{m=n} et les troisi\`eme et cinqui\`eme si \m{l=m}. Ces applications
jouent le r\^ole de la composition des morphismes dans le cas des
faisceaux. On suppose qu'elles sont toutes surjectives et qu'elles
v\'erifient les propri\'et\'es usuelles qu'on attend des applications
habituelles de composition. Cela signifie que les diagrammes suivants
sont commutatifs (les fl\`eches \'etant les fl\`eches \'evidentes) :
$$\diagram{
A_{kj}\ot A_{ji}\ot A_{ih} & \hfl{}{} & A_{ki}\ot A_{ih} \cr
\vfl{}{} & & \vfl{}{}\cr
A_{kj}\ot A_{jh} & \hfl{}{} & A_{kh} \cr
}\ \ \ \
\diagram{
H_{lk}\ot A_{kj}\ot A_{ji} & \hfl{}{} & H_{lj}\ot A_{ji} \cr
\vfl{}{} & & \vfl{}{}\cr
H_{lk}\ot A_{ki} & \hfl{}{} & H_{li}\cr
}$$
$$\diagram{
B_{ml}\ot H_{lj}\ot A_{ji} & \hfl{}{} & H_{mj}\ot A_{ji} \cr
\vfl{}{} & & \vfl{}{}\cr
B_{ml}\ot H_{li} & \hfl{}{} & H_{mi}\cr
}\ \ \ \
\diagram{
B_{nm}\ot B_{ml}\ot H_{li} & \hfl{}{} & B_{nl}\ot H_{li} \cr
\vfl{}{} & & \vfl{}{}\cr
B_{nm}\ot H_{mi} & \hfl{}{} & H_{ni}\cr
}$$
$$\diagram{
B_{on}\ot B_{nm}\ot B_{ml} & \hfl{}{} & B_{om}\ot B_{ml} \cr
\vfl{}{} & & \vfl{}{}\cr
B_{on}\ot B_{nl} & \hfl{}{} & B_{ol}\cr
}$$
On supposera aussi que les applications
$$H_{li}^* \ot A_{ji}\lra H_{lj}^* \ \ , \ \
H_{mi}^*\ot B_{ml}\lra H_{li}^*$$
induites par les applications de composition sont surjectives. Soient
\m{M_i, 1\leq i\leq r,}
\break\m{N_l, 1\leq l\leq s} des espaces vectoriels de
dimension finie. On notera
$$m_i = \dim(M_i), \ n_l = \dim(N_l), \ 1\leq i\leq r, 1\leq l\leq s.$$
On veut \'etudier l'espace vectoriel
$$W = \bigoplus_{1\leq i\leq r,1\leq l\leq s}
\Hom(H_{li}^*\ot M_i,N_l) .$$

\subsection{D\'efinition du groupe \m{G}}

Soit \m{G_L} l'ensemble des matrices
$$g = \pmatrix{g_1 & 0   & .       & . & . & 0   \cr
u_{21}             & g_2 & .       & . & . & 0   \cr
.                  & .   & .       &   &   & .   \cr
.                  &     &         & . &   & .   \cr
.                  &     & u_{ij}  &   & . & .   \cr
u_{r1}             & .   & .       & . & . & g_r \cr
} ,$$
avec \ \m{g_i\in GL(M_i)}, et pour \ \m{1\leq j < i \leq r},
$$u_{ij}\in \Hom(A_{ij}^*\ot M_j,M_i) = \Hom(M_j,A_{ij}\ot M_i) .$$
Soit \m{G_R} l'ensemble des matrices
$$h = \pmatrix{h_1 & 0   & .       & . & . & 0   \cr
v_{21}             & h_2 & .       & . & . & 0   \cr
.                  & .   & .       &   &   & .   \cr
.                  &     &         & . &   & .   \cr
.                  &     & v_{lm}  &   & . & .   \cr
v_{s1}             & .   & .       & . & . & h_s \cr
} ,$$
avec \ \m{h_l\in GL(N_l)}, et pour \ \m{1\leq m < l \leq s},
$$v_{lm}\in \Hom(B_{lm}^*\ot N_m,N_l) = \Hom(N_m,B_{lm}\ot N_l) .$$
On d\'efinit une loi de composition, not\'ee \m{*}, de la fa\c con
suivante : si
$$ u_{kj}\in \Hom(A_{kj}^*\ot M_j,M_k) \
\ \ {\rm et} \ \ u_{ji}\in \Hom(A_{ji}^*\ot M_i,M_j),$$
alors
$$u_{kj}*u_{ji}\in L(A_{ki}^*\ot M_i,M_k) $$
est la composition
$$A_{ki}^*\ot M_i\hfl{u_{ji}}{} A_{ki}^*
\ot A_{ji}\ot M_j\hfl{}{}A_{kj}^*\ot M_j\hfl{u_{kj}}{}M_k ,$$
ou l'application du milieu est induite par la composition
$$A_{kj}\ot A_{ji}\lra A_{ki}.$$
On d\'efinit une structure de groupe sur \m{G_L} de la fa\c con suivante :
si \m{g,g' \in G_L}, avec
$$g = \pmatrix{g_1 & 0   & .       & . & . & 0   \cr
u_{21}             & g_2 & .       & . & . & 0   \cr
.                  & .   & .       &   &   & .   \cr
.                  &     &         & . &   & .   \cr
.                  &     & u_{ij}  &   & . & .   \cr
u_{r1}             & .   & .       & . & . & g_r \cr
} \ , \
g' = \pmatrix{g'_1 & 0   & .       & . & . & 0   \cr
u'_{21}             & g'_2 & .       & . & . & 0   \cr
.                  & .   & .       &   &   & .   \cr
.                  &     &         & . &   & .   \cr
.                  &     & u'_{ij}  &   & . & .   \cr
u'_{r1}             & .   & .       & . & . & g'_r \cr
} ,$$
alors
$$g'g = \pmatrix{g''_1 & 0   & .       & . & . & 0   \cr
u''_{21}             & g''_2 & .       & . & . & 0   \cr
.                  & .   & .       &   &   & .   \cr
.                  &     &         & . &   & .   \cr
.                  &     & u''_{ij}  &   & . & .   \cr
u''_{r1}             & .   & .       & . & . & g''_r \cr
} ,$$
avec
$$g''_i = g'_i\circ g_i \ \ (1\leq i\leq r),$$
$$u''_{ij} = u'_{ij}\circ g_j + \sum_{1\leq k<i-j} u'_{i,j+k}*u_{j+k,j}
+ g'_i\circ u_{ji}  \ \ (1\leq j < i \leq r) .$$
La v\'erification qu'on obtient ainsi une structure de groupe sur
\m{G_L} est imm\'ediate. On d\'efinit une structure de groupe analogue
sur \m{G_R}. Soit
$$G = G_L\times G_R .$$

\subsection{D\'efinition de l'action de \m{G} sur \m{W}}

On va d\'efinir une action \`a gauche de \m{G_L} sur \m{W} et une action
\`a droite de \m{G_R} sur \m{W}. L'action de \m{G} sur \m{W} en d\'ecoule :
si \m{(g,h)\in G} et \m{w\in W}, on a
$$(g,h).w \ = \ h.w.g^{-1} .$$
Soit \ \m{w = (\phi_{li})_{1\leq i\leq r,1\leq l\leq s}\in W} (donc
\m{\phi_{il}} est un application lin\'eaire
 \ \m{H_{li}^*\ot M_i\lra N_l}). Soit
$$g = \pmatrix{g_1 & 0   & .       & . & . & 0   \cr
u_{21}             & g_2 & .       & . & . & 0   \cr
.                  & .   & .       &   &   & .   \cr
.                  &     &         & . &   & .   \cr
.                  &     & u_{ij}  &   & . & .   \cr
u_{r1}             & .   & .       & . & . & g_r \cr
} $$
un \'el\'ement de \m{G_L}. Alors \m{w.g =
(\phi'_{li})_{1\leq i\leq r,1\leq l\leq s}}, o\`u
$$\phi'_{li} = \sigg_{i\leq j\leq r}\psi_{ijl},$$
 \m{\psi_{iil}} \'etant la composition
$$M_i\ot H_{li}^*\hfl{g_i}{}M_i\ot H_{li}^*\hfl{\phi_{li}}{}N_l$$
et, si \m{i<j\leq r}, \m{\psi_{ijl}} la composition
$$M_i\ot H_{li}^*\hfl{u_{ji}}{}M_j\ot A_{ji}\ot H_{li}^*
\hfl{}{}M_j\ot H_{lj}^*\hfl{\phi_{lj}}{}N_l,$$
l'application du milieu  \'etant induite par la composition \
\m{H_{lj}\ot A_{ji}\lra H_{li}}.
L'action de \m{G_R} est analogue.

\subsection{Notions de (semi-)stabilit\'e}

On veut d\'efinir une notion de {\em (semi-)stabilit\'e} pour les points
de \m{W}. On ne peut pas appliquer la g\'eom\'etrie invariante si
\m{r>1} ou \m{s>1} car le groupe \m{G} n'est pas r\'eductif. On va d\'efinir
deux sous-groupes canoniques de \m{G}.
Soit \m{H_L} (resp. \m{G_{L,red}} ) le sous-groupe de \m{G_L} form\'e des
\'el\'ements
$$\pmatrix{g_1 & 0   & .       & . & . & 0   \cr
u_{21}             & g_2 & .       & . & . & 0   \cr
.                  & .   & .       &   &   & .   \cr
.                  &     &         & . &   & .   \cr
.                  &     & u_{ij}  &   & . & .   \cr
u_{r1}             & .   & .       & . & . & g_r \cr
} $$
tels que \ \m{g_i = I_{M_i}} \ pour \m{1\leq i\leq r} (resp. \
\m{u_{ij} = 0} \ pour \m{1\leq j < i\leq r}).
Alors \m{H_L} est un sous-groupe unipotent normal maximal
de \m{G_L}, \m{G_{L,red}} est un sous-groupe r\'eductif de \m{G_L}
et l'inclusion \m{G_{L,red}\subset G_L} induit un isomorphisme
\m{G_{L,red}\simeq G_L/H_L}. On d\'efinit de m\^eme les sous-groupes \m{H_R}
et \m{G_{R,red}} de \m{G_R}.

Maintenant soient
$$H = H_L\times H_R \ , \ G_{red} = G_{L,red}\times G_{R,red} .$$
Alors \m{H} est un sous-groupe unipotent normal maximal de \m{G} et
\m{G_{red}} est un sous-groupe r\'eductif de \m{G}.

L'action de \m{G_{red}} sur \m{W} est un cas particulier des actions
trait\'ees dans \cite{king}. Soient
\m{\lambda_1,\ldots,\lambda_r,} \m{\mu_1,\ldots,\mu_s} des nombres
rationnels positifs tels que
$$\sigg_{1\leq i\leq r}\lambda_im_i = \sigg_{1\leq l\leq s}\mu_ln_l.$$

\begin{defin}
On dit qu'un \'el\'ement \m{(\phi_{li})} de \m{W} est
\m{G_{red}}-semi-stable ( resp. \m{G_{red}}-stable) relativement \`a
\m{(\lambda_1,\ldots,\lambda_r,} \m{\mu_1,\ldots,\mu_s)} si la propri\'et\'e
suivante est v\'erifi\'ee : soient \m{M'_i\subset M_i},
\m{N'_l\subset N_l} des sous-espaces vectoriels tels
que l'un au moins des \m{N'_l} soit distinct de \m{N_l} et que
pour \m{1\leq i\leq r}, \m{1\leq l\leq s}, on ait
$$\phi_{li}(H_{li}^*\ot M'_i)\subset N'_l.$$
Alors on a
$$\sigg_{1\leq i\leq r}\lambda_i\dim(M'_i)\ \leq \
\sigg_{1\leq l\leq s}\mu_l\dim(N'_l) \ \ \ {\rm (resp. \ } <{\rm )} \ .$$
\end{defin}

\begin{defin}
On dit qu'un \'el\'ement \m{x} de \m{W} est
\m{G}-semi-stable ( resp. \m{G}-stable) relativement \`a
\m{(\lambda_1,\ldots,\lambda_r,} \m{\mu_1,\ldots,\mu_s)} si tous les
points de l'orbite \m{H.x} sont \m{G_{red}}-semi-stables ( resp.
\m{G_{red}}-stables) relativement \`a
\m{(\lambda_1,\ldots,\lambda_r,} \m{\mu_1,\ldots,\mu_s)}.
\end{defin}

On note \m{W^{ss}(\lambda_1,\ldots,\lambda_r,\mu_1,\ldots,\mu_s)}
(resp. \m{W^{s}(\lambda_1,\ldots,\lambda_r,\mu_1,\ldots,\mu_s)}), ou
plus simplement \m{W^{ss}} (resp. \m{W^{s}})
si aucune confusion n'est \`a craindre, l'ouvert
de \m{W} constitu\'e des points
\m{G}-semi-stables ( resp. \m{G}-stables) relativement \`a
\m{(\lambda_1,\ldots,\lambda_r,} \m{\mu_1,\ldots,\mu_s)}.

\subsection{Cas d'existence d'un bon quotient projectif}

On donne dans \cite{dr_tr} des conditions suffisantes portant sur
\m{\lambda_1,\ldots,\lambda_r,} \m{\mu_1,\ldots,\mu_s}, pour qu'il existe
un bon quotient
$$\pi : W^{ss}\lra M = M(\lambda_1,\ldots,\lambda_r,\mu_1,\ldots,\mu_s)$$
par \m{G} avec \m{M} projective. Dans ce cas \m{M} est normale et la
restriction de \m{\pi}
$$W^{s}\lra M^s = \pi(W^s)$$
est un quotient g\'eome\'trique. Le r\'esultat le plus g\'en\'eral est assez
compliqu\'e. Rappelons simplement ici le cas des morphismes de type
\m{(2,1)}, le seul qu'on utilisera ici (dans le \para 8).
Il faut d'abord d\'efinir certaines constantes. Soit \m{k>0}
un entier. Soient
$$\tau : H_{11}^*\ot A_{21}\lra H_{12}^*$$
l'application lin\'eaire d\'eduite de la composition \
\m{H_{12}\ot A_{21}\lra H_{11}}, et
$$\tau_k = \tau_1\ot I_{\scx{k}} :
H_{11}^*\ot(A_{21}\ot\cx{k})\lra H_{12}^*\ot\cx{k}.$$
Soit \m{\cal K} l'ensemble des sous-espaces vectoriels propres
\ \m{K\subset A_{21}\ot\cx{k}} \
tels que pour tout sous-espace propre \m{F\subset\cx{k}}, \m{K} ne soit pas
contenu dans  \m{A_{21}\ot F}. Alors posons
$$c(\tau, k) = \supp_{K\in{\cal K}}(\q{\codim(\tau_k(H_{11}^*\ot K)}
{\codim(K)}).$$

Dans le cas des morphismes de type \m{(2,1)}, les notions de
semi-stabilit\'e sont d\'efinies \`a partir de triplets
$$(\lambda_1,\lambda_2,\q{1}{n_1})$$
tels que
\ \m{\lambda_1 m_1 + \lambda_2 m_2 = 1}.
Elles d\'ependent donc essentiellement d'un param\`etre. Le r\'esultat
suivant est d\'emontr\'e dans \cite{dr_tr} :

\begin{xtheo}
Il existe un bon quotient projectif
\m{W^{ss}//G} d\`es que
$$\q{\lambda_2}{\lambda_1}>\dim(A_{21}) \ \ \ {\rm et} \ \
\lambda_2\geq \q{\dim(A_{21})}{n_1} c(\tau,m_2).$$
\end{xtheo}

\subsection{Dualit\'e}

La notion de dualit\'e est claire dans le contexte des morphismes de
faisceaux (en supposant qu'ils sont localement libres). Au lieu
d'\'etudier des morphismes
\ \m{{\cal E}\lra{\cal F}} \
on consid\`ere les morphismes transpos\'es
\ \m{{\cal F}^*\lra{\cal E}^*}.
Dans le cas g\'en\'eral, on pose \m{r'=s, s'=r},

\noindent\m{A'_{ij}=B_{s+1-j,s+1-i}},
\m{B'_{lm}=A_{r+1-m,r+1-l}},

\noindent\m{H'_{li}=H_{s+1-i,r+1-l}}, les compositions
sont les m\^emes. On prend \m{M'_i = N_{s+1-i}^*},

\noindent\m{N'_l=M_{r+1-l}^*}.
L'espace associ\'e \m{W'} est isomorphe \`a \m{W}, le facteur
\m{\Hom(M_i\ot H_{li},N_l)} s'identifiant \`a \m{\Hom(M'_{s+1-l}\ot
H'_{r+1-i,s+1-l},N'_{r+1-i})}. Le groupe \m{G'} est le m\^eme (sauf pour
l'ordre des facteurs, c'est-\`a-dire \m{G'_L = G_R} et \m{G'_R = G_L}). Les
actions des groupes sont bien s\^ur les m\^emes.
\bigskip
\bigskip

\section{Mutations d\'efinies \`a l'aide de la suite spectrale de
Beilinson g\'en\'eralis\'ee}

\subsection{Rappels sur les suites spectrales de Beilinson
g\'en\'eralis\'ees sur les espaces projectifs}

Les d\'efinitions et propri\'et\'es de base des h\'elices de \fexcs sur
\m{\proj{n}} se trouvent dans \cite{go_ru}.

\subsubsection{H\'elices de \fexcs sur \m{\proj{n}}}

Une {\em h\'elice} $\gamma = (E_i)_{i \in Z \hskip -4pt Z}$ de \fexcs
sur \m{\proj{n}} poss\`ede les propri\'et\'es suivantes :

\medskip
\noindent 1) C'est une suite {\em p\'eriodique}, c'est-\`a-dire qu'on a
\ \m{E_{i+n+1} \simeq E_i(n+1)} \ pour tout entier \m{i}.

\medskip
\noindent 2) On a \ \m{\chi(E_i,E_j) = 0} \ si \ \m{j < i \leq j+n} .

\medskip
\noindent 3) Pour tout entier \m{i}, le morphisme canonique
$$ev : E_{i-1}\otimes \Hom(E_{i-1},E_i)\lra E_i, \ \ \
\rm{( resp. \ } ev^* : E_i\lra E_{i+1}\otimes\Hom(E_i,E_{i+1})^* \rm{ )}$$
\noindent est surjectif (resp. injectif) et son noyau (resp. conoyau) est un
\fexc $E$ (resp. $F$). De plus, la suite p\'eriodique de fibr\'es vectoriels
bas\'ee sur
$$(E_{i-2},E,E_{i-1},E_{i+1},\ldots,E_{i+n-2})
\ \ \ {\rm (resp. } \ \
(E_{i-1},E_{i+1},F,E_{i+2},\ldots,E_{i+n-1})\ {\rm)}$$
 est une h\'elice
(le terme $E$ (resp. $F$) \'etant d'indice \m{i-1} (resp. \m{i+1})).
Cette h\'elice s'appelle {\em mutation \`a gauche} (resp. {\em mutation \`a
droite} de $\gamma$ en \m{E_i}, et est not\'ee \m{L_{E_i}(\gamma)} (resp.
\m{R_{E_i}(\gamma)}). Le \fexc \m{E} (resp. \m{F}) est not\'e
\m{L_\gamma(E_i)} (resp. \m{R_\gamma(E_i)}).

\medskip
\noindent 4) On pose \ \m{L^2_{E_i}(\gamma) = L_E\circ L_{E_i}(\gamma)}.
C'est une h\'elice ayant pour base une suite de la forme
$$(E_{i-3},E',E_{i-2},E_{i-1},E_{i+1},\ldots,E_{i+n-3}).$$
 On d\'efinit de
m\^eme \m{L^p_{E_i}(\gamma)} pour tout entier p tel que \
\m{1\leq p < n}. C'est la suite infinie p\'eriodique bas\'ee sur une
suite du type
$$(E_{i-p-1},E^{(p)},E_{i-p},\ldots,E_{i-1},E_{i+1},\ldots,
E_{i+n-p-1}),$$
 \m{E^{(p)}} \'etant un \fexc et d'indice \m{i-p}. En
particulier, \m{L^{n-1}_{E_i}(\gamma)} est bas\'ee sur la suite
$$(E_{i-n},E^{(n-1)},E_{i-n+1},\ldots,E_{i-1}). $$
L'h\'elice
$$L^n_{E_i}(\gamma) = L_{E^{(n-1)}}\circ L^{(n-1)}(\gamma)$$
est bas\'ee sur une suite du type
$$(E^{(n)},E_{i-n},\ldots,E_{i-1}),$$
\m{E^{(n)}} \'etant un \fexc, d'indice \m{i-n}. Alors on a
$$E^{(n)} \simeq E_{i-n-1},$$
c'est-\`a-dire que \m{L^n_{E_i}(\gamma)} est \'egale \`a \m{\gamma} \`a un
d\'ecalage pr\`es. On notera
$$L^p_\gamma(E_i) = E^{(p)}, $$
et en consid\'erant les mutations \`a droite on d\'efinit de m\^eme les
\fexcs \m{R^p_\gamma(E_i)}.

\medskip
\noindent 5) On a \ $L_{E_i}\circ L^2_{E_{i+1}}(\gamma) = L^2_{E_{i+1}}\circ
L_{E_i}(\gamma)$ .

\bigskip

On a bien s\^ur des propri\'et\'es analogues \`a 4) et 5) concernant les
mutations \`a droite. L'h\'elice la plus simple est
$$({\cal O}(i))_{i \in Z \hskip -4pt Z}$$
et toutes les h\'elices de \fexcs connues
peuvent s'obtenir en partant de cette h\'elice et en lui faisant subir une
suite finie de mutations et une translation des indices.

Une {\em base} de l'h\'elice
\m{\gamma = (E_i)_{i \in Z \hskip -4pt Z}} de \fexcs sur \m{\proj{n}} est
une suite
$$\sigma = (E_i,\ldots,E_{i+n})$$
extraite de \m{\gamma}. A cause de la propri\'et\'e 1-, \m{\gamma} peut
\^etre reconstitu\'ee \`a partir de \m{\sigma}. Les
notions de mutations \`a droite et \`a gauche s'\'etendent de mani\`ere
\'evidente aux bases d'h\'elice. Si \m{i< j< n}, on note
$$L_{j+1}^1(\sigma) =
(E_i,\ldots, E_{j-1}, L^1_\gamma(E_{j+1}), E_j,E_{j+2},\ldots, E_n),$$
$$R_j^1(\sigma) =
(E_i,\ldots, E_{j-1}, E_{j+1}, R^1_\gamma(E_j),E_{j+2},\ldots, E_n).$$
Plus g\'en\'eralement, si \m{1\leq p\leq j-i}, on pose
$$L_{j+1}^p(\sigma) =
(E_i,\ldots, E_{j-p}, L^p_\gamma(E_{j+1}),E_{j-p+1},\ldots,
E_j,E_{j+2},\ldots, E_n),$$
et si \m{1\leq q\leq n-j-1},
$$R_j^q(\sigma) = (E_i,\ldots, E_{j-1}, E_{j+1},\ldots, E_{j+q},
R^q_\gamma(E_j),E_{j+q+1},\ldots, E_n).$$

\subsubsection{Suite spectrale de Beilinson g\'en\'eralis\'ee}

\subsubsubsection{4.1.2.1 -}{D\'efinition}

Soit \ \m{\sigma = (E_0,\ldots,E_n)} \ une base d'h\'elice sur \m{\proj{n}}.
On associe \`a  \m{\sigma} une autre base d'h\'elice, dite {\em duale} de
\m{\sigma}, et not\'ee
$$\sigma^* = (E_{\sigma 0},\ldots,E_{\sigma n}),$$
d\'efinie par
$$E_{\sigma p} = L^p_\sigma(E_{p})^*(-n-1) = R^{n-p}_\sigma(E_{p})^* .$$
C'est une base d'une autre h\'elice que \m{\gamma}.
Si \m{\gamma} est l'h\'elice engendr\'ee par \m{\sigma}, on note
\m{\gamma^*} l'h\'elice engendr\'ee par \m{\sigma^*}.
On montre qu'il existe
une r\'esolution canonique de la diagonale
\m{\Delta} de \ \m{\proj{n}\times\proj{n}} :
$$0\lra E_0\timex E_{\sigma 0}\lra\cdots\lra E_{n}\timex E_{\sigma,n}
\lra{\cal O}_{\Delta}\lra 0 .$$
On en d\'eduit, pour tout faisceau coh\'erent \m{\E} sur \m{\proj{n}}, une
suite spectrale \m{E_r^{pq}} de faisceaux coh\'erents sur \m{\proj{n}},
convergeant vers \m{\E} en degr\'e \m{0} et vers \m{0} en tout autre
degr\'e, et dont les termes \m{E_1^{p,q}} pouvant \'eventuellement \^etre
non nuls sont les
$$E_1^{p,q} = E_{p+n}\otimes H^q(E_{\sigma,p+n}\otimes\E) \ ,
\ -n \leq p\leq 0 \ , \ 0\leq q\leq n .$$
On en d\'eduit le {\em complexe de Beilinson} :
$$0\lra\F_{-n}\lra\F_{-n+1}\lra\ldots\lra\F_{n-1}\lra
\F_n\lra 0,$$
o\`u \ \m{\F_i = \som_{p+q=i} E_1^{pq}}. Il est exact en degr\'e
diff\'erent
de \m{0}, et sa cohomologie en degr\'e \m{0} est isomorphe \`a \m{\E}.

\bigskip

\subsubsubsection{4.1.2.2 -}{Bases duales et mutations}

Soit \ \m{\sigma = (E_0,\ldots,E_n)} \ une base d'h\'elice, et \m{j} un
entier tel que \ \m{1\leq j< n}. Alors on a
$$L_{j+1}(\sigma)^* = R_j(\sigma^*),
\ \ \ R_j(\sigma)^* = L_{j+1}(\sigma^*).$$

\bigskip

\subsection{Mutations de morphismes}

Soient \m{r,s} des entiers positifs, et \m{n=r+s-1}.
Soit
$$\sigma = (E_1,\ldots,E_r,F_1,\ldots,F_s)$$
une base d'h\'elice de fibr\'es exceptionnels sur
$\proj{n}$, \m{\gamma} l'h\'elice engendr\'ee par \m{\sigma}.
Pour \m{1\leq i\leq r}, le morphisme canonique de fibr\'es vectoriels
$$E_i\lra F_1\ot \Hom(E_i,F_1)^*$$
est surjectif, et son conoyau \m{G_i} est un fibr\'e exceptionnel.
On va d\'efinir une suite \m{\sigma_r}, \m{\sigma_{r-1}}, \m{\ldots},
\m{\sigma_0} de bases d'h\'elice par
$$\sigma_r=\sigma,$$
et si \m{p} est un entier tel que \m{0\leq p\leq r-1},
$$\sigma_p = R_{p}(\sigma_{p+1}).$$
On a
$$\sigma_p =
(E_1,\ldots,E_p,F_1,G_{p+1},\ldots,G_r,F_2,\ldots, F_s).$$
Il d\'ecoule du \para 3.2.2 qu'on a
$$\sigma_p^* = (E_{\sigma 0},\ldots,E_{\sigma,p-1},
L^p_{\gamma^*}(E_{\sigma r}),E_{\sigma p},\ldots, E_{\sigma n}).$$

En utilisant la suite spectrale de Beilinsion g\'en\'eralis\'ee associ\'ee
\`a \m{\sigma}, on d\'emontre ais\'ement ce qui suit : soit ${\cal U}$ un
faisceau coh\'erent sur $\proj{n}$ tel que
$$H^j({\cal U}\ot E_{\sigma i}) = \lbrace 0\rbrace$$
si \m{0\leq i<r} et \m{j\not = n-i-1}, ou \m{r\leq i\leq n} et
\m{j\not = n-i}. On pose
$$M_i = H^{n-i}({\cal U}\ot E_{\sigma,i-1}) \ \ \ {\rm pour\ \ }
1\leq i\leq r,$$
$$N_l = H^{n-r-l+1}({\cal U}\ot E_{\sigma,l+r-1}) \ \ \ {\rm pour\ \ }
1\leq l\leq s,$$
de telle sorte que le diagramme de Beilinson de \m{\cal U} a l'allure
suivante
$$\matrix{
0        & . & . & . & 0         & 0        & . & . & . & 0        \cr
M_1      &   &   &   & 0         & 0        & . & . & . & 0        \cr
0        & . &   &   & .         & .        &   &   &   & .        \cr
.        &   & . &   & .         & .        &   &   &   & .        \cr
.        &   &   & . & 0         & 0        &   &   &   & .        \cr
0        & . & . & 0 & M_r       & N_1      & 0 & . & . & 0        \cr
.        &   &   &   & 0         & 0        & . &   &   & .        \cr
.        &   &   &   & .         & .        &   & . &   & .        \cr
.        &   &   &   & .         & .        &   &   & . & .        \cr
0        &   &   &   & .         & .        &   &   &   & N_s      \cr
}$$
Alors il existe une suite exacte
$$0\lra\som_{1\leq i\leq r}(E_i\ot M_i) \hfl{\Phi}{}
\som_{1\leq l\leq s}(F_l\ot N_l)\lra{\cal U}\lra 0.$$
Soit \m{p} un entier tel que \m{0\leq p\leq r-1}. On note
$$f_p : \som_{p+1\leq j\leq r}(\Hom(E_j,F_1)^*\ot M_j)\lra N_1$$
l' application lin\'eaire d\'eduite de \m{\Phi}.
Alors, si $f_p$ est surjective, on peut montrer que la suite spectrale
de Beilinson g\'en\'eralis\'ee associ\'ee a \m{\sigma_p},
appliqu\'ee \`a ${\cal U}$, donne une suite exacte
$$0\lra \biggl(\bigoplus_{1\leq i\leq p}(E_i\ot M_i)\biggr)\oplus
(F_1\ot\ker(f_p))\lra
\biggr(\bigoplus_{p<j\leq r}(G_j\ot M_j)\biggl)\oplus
\biggr(\bigoplus_{2\leq l\leq s}(F_l\ot N_l)\biggl)\lra{\cal U}\lra 0.$$
On a aussi bien s\^ur un \'enonc\'e r\'eciproque.
Ce r\'esultat va \^etre g\'en\'eralis\'e dans le chapitre suivant.

\bigskip
\bigskip

\section{Mutations en termes de morphismes de faisceaux}

On va d'abord d\'emontrer deux r\'esultats, qu'on appliquera ensuite \`a
la d\'efinition des mutations de morphismes de type \m{(r,s)}. On \'etudie
des faisceaux coh\'erents pouvant \^etre repr\'esent\'es comme conoyaux
de morphismes injectifs de faisceaux d'un certain type. Une \'etude
similaire pourrait sans doute \^etre faite sur les noyaux.

\subsection{R\'esultats g\'en\'eraux}

Soient \m{\E}, \m{\E'}, \m{\F}, \m{\F'} et \m{\Gamma} des faisceaux
coh\'erents sur une vari\'et\'e projective irr\'eductible \m{X}, avec
\m{\Gamma} simple. On
suppose que le morphisme canonique
$$ev : \Gamma\ot \Hom(\Gamma,\F)\lra \F$$
est surjectif. Soit \m{\E_0} son noyau. On suppose que le morphisme
canonique
$$ev^* : \E'\lra \Gamma\ot\Hom(\E',\Gamma)^*$$
est injectif. Soit \m{\F_0} son conoyau. On suppose enfin que
$$\Hom(\E',\E_0) = \Ext^1(\E',\E_0) = \Ext^1(\F_0,\F') =
\Ext^1(\E,\E_0) = \lbrace 0\rbrace.$$

De la suite exacte
$$0\lra\E_0\lra\Gamma\ot\Hom(\Gamma,\F)\lra\F\lra 0$$
on d\'eduit un isomorphisme
$$\Hom(\E',\F) \ \simeq \ \Hom(\Hom(\E',\Gamma)^*,\Hom(\Gamma,\F)).$$
Si \ \m{\lambda\in\Hom(\Hom(\E',\Gamma)^*,\Hom(\Gamma,\F))}, le morphisme
\ \m{\E'\lra\F} \ correspondant est la compos\'ee
$$\E'\hfl{ev^*}{}\Gamma\ot\Hom(\E',\Gamma)^*\hfl{I_\Gamma\ot\lambda}{}
\Gamma\ot\Hom(\Gamma,\F)\hfl{ev}{}\F.$$

\bigskip

\begin{xprop}
Soit
$$\Phi : \E\oplus\E'\lra\F\oplus\F'$$
un morphisme injectif de faisceaux, et \m{\cal U} son conoyau. Soit
$$\lambda : \Hom(\E',\Gamma)^*\lra\Hom(\Gamma,\F)$$
l'application lin\'eaire d\'eduite du morphisme \ \m{\E'\lra\F} \ d\'efini
par \m{\Phi}.

\medskip

\noindent 1 - On suppose que \m{\lambda} est surjective et que
$$\Ext^1(\E,\Gamma) = \lbrace 0\rbrace.$$
Alors il existe une suite exacte
$$0\lra \E\oplus\E_0\oplus(\Gamma\ot\ker(\lambda))\lra\F_0\oplus\F'
\lra{\cal U}\lra 0.$$

\medskip

\noindent 2 - On suppose que \m{\lambda} est injective et que
$$\Ext^1(\Gamma,\F_0) = \Ext^1(\Gamma,\F') = \lbrace 0\rbrace.$$
Alors il existe une suite exacte
$$0\lra \E\oplus\E_0\lra(\Gamma\ot\coker(\lambda))\oplus\F_0\oplus\F'
\lra{\cal U}\lra 0.$$
\end{xprop}

\noindent{\em D\'emonstration}. On consid\`ere le morphisme
$$A : \E'\lra (\Gamma\ot\Hom(\E',\Gamma)^*)\oplus\F' = {\cal A}$$
dont la premi\`ere composante est \m{ev^*} et la seconde provient de
\m{\Phi}. \Oad

\bigskip

\begin{picture}(360,230)
\put(135,220){$0$}
\put(280,220){$0$}
\put(10,170){$0$}
\put(55,170){$\E'$}
\put(115,170){$\Gamma\ot\Hom(\E',\Gamma)^*$}
\put(280,170){$\F_0$}
\put(365,170){$0$}
\put(10,120){$0$}
\put(55,120){$\E'$}
\put(135,120){$\cal A$}
\put(270,120){$\coker(A)$}
\put(365,120){$0$}
\put(135,70){$\F'$}
\put(280,70){$\F'$}
\put(135,20){$0$}
\put(280,20){$0$}

\put(137,184){\vector(0,1){30}}
\put(137,134){\vector(0,1){30}}
\put(137,84){\vector(0,1){30}}
\put(137,34){\vector(0,1){25}}
\put(282,184){\vector(0,1){30}}
\put(282,134){\vector(0,1){30}}
\put(282,84){\vector(0,1){30}}
\put(282,34){\vector(0,1){30}}
\put(56,163){\line(0,-1){29}}
\put(58,163){\line(0,-1){29}}

\put(70,173){\vector(1,0){35}}
\put(295,173){\vector(1,0){63}}
\put(23,173){\vector(1,0){28}}
\put(204,173){\vector(1,0){67}}
\put(23,123){\vector(1,0){28}}
\put(70,123){\vector(1,0){56}}
\put(90,126){A}
\put(152,123){\vector(1,0){111}}
\put(323,123){\vector(1,0){35}}
\put(150,73){\line(1,0){122}}
\put(150,71){\line(1,0){122}}
\end{picture}

\bigskip

Puisque \ \m{\Ext^1(\F_0,\F') = \lbrace 0\rbrace}, on a un isomorphisme
$$\coker(A) \ \simeq \ \F'\oplus\F_0 .$$

On suppose maintenant que les hypoth\`eses de 1- sont v\'erifi\'ees. Soit
$$\pi : \Gamma\ot\Hom(\E',\Gamma)^*\lra\F$$
le morphisme compos\'e
$$\Gamma\ot\Hom(\E',\Gamma)^*\hfl{I_\Gamma\ot\lambda}{}
\Gamma\ot\Hom(\Gamma,\F)\hfl{ev}{}\F.$$
Alors on a
$$\ker(\pi)\simeq\E_0\oplus(\ker(\lambda)\ot\Gamma).$$
Soit \m{\cal V} le conoyau du morphisme injectif
$$\E'\lra\F\oplus\F'$$
d\'eduit de \m{\Phi}. \Oad

\bigskip

\begin{picture}(360,230)
\put(135,220){$0$}
\put(280,220){$0$}
\put(10,170){$0$}
\put(55,170){$\E'$}
\put(123,170){$\F\oplus\F'$}
\put(280,170){$\cal V$}
\put(365,170){$0$}
\put(10,120){$0$}
\put(55,120){$\E'$}
\put(135,120){$\cal A$}
\put(262,120){$\F'\oplus\F_0$}
\put(365,120){$0$}
\put(125,70){$\ker(\pi)$}
\put(270,70){$\ker(\pi)$}
\put(135,20){$0$}
\put(280,20){$0$}

\put(137,184){\vector(0,1){30}}
\put(137,134){\vector(0,1){30}}
\put(137,84){\vector(0,1){30}}
\put(137,34){\vector(0,1){25}}
\put(140,149){$\pi\oplus I_{\F'}$}
\put(282,184){\vector(0,1){30}}
\put(282,134){\vector(0,1){30}}
\put(282,84){\vector(0,1){30}}
\put(282,34){\vector(0,1){30}}
\put(56,163){\line(0,-1){29}}
\put(58,163){\line(0,-1){29}}

\put(70,173){\vector(1,0){45}}
\put(295,173){\vector(1,0){63}}
\put(23,173){\vector(1,0){28}}
\put(170,173){\vector(1,0){101}}
\put(23,123){\vector(1,0){28}}
\put(70,123){\vector(1,0){56}}
\put(90,126){A}
\put(152,123){\vector(1,0){103}}
\put(312,123){\vector(1,0){46}}
\put(160,73){\line(1,0){102}}
\put(160,71){\line(1,0){102}}
\end{picture}

\bigskip

On a une suite exacte
$$0\lra\E\lra{\cal V}\lra{\cal U}\lra 0,$$
et l'inclusion \ \m{\E\subset{\cal V}} \ se rel\`eve en un morphisme
injectif
$$\E\lra\F'\oplus\F_0$$
(car \ \m{\Ext^1(\E,\E_0) = \Ext^1(\E,\Gamma) = \lbrace 0\rbrace}). On
note \m{\cal W} le conoyau de ce morphisme. On
a alors un diagramme commutatif avec lignes et colonnes exactes, dont
la ligne verticale du milieu provient du diagramme pr\'ec\'edent :

\bigskip

\begin{picture}(360,230)
\put(135,220){$0$}
\put(280,220){$0$}
\put(10,170){$0$}
\put(55,170){$\E$}
\put(133,170){$\cal V$}
\put(280,170){$\cal U$}
\put(365,170){$0$}
\put(10,120){$0$}
\put(55,120){$\E$}
\put(120,120){$\F'\oplus\F_0$}
\put(280,120){$\cal W$}
\put(365,120){$0$}
\put(125,70){$\ker(\pi)$}
\put(270,70){$\ker(\pi)$}
\put(135,20){$0$}
\put(280,20){$0$}

\put(137,184){\vector(0,1){30}}
\put(137,134){\vector(0,1){30}}
\put(137,84){\vector(0,1){30}}
\put(137,34){\vector(0,1){25}}
\put(282,184){\vector(0,1){30}}
\put(282,134){\vector(0,1){30}}
\put(282,84){\vector(0,1){30}}
\put(282,34){\vector(0,1){30}}
\put(56,163){\line(0,-1){29}}
\put(58,163){\line(0,-1){29}}

\put(70,173){\vector(1,0){55}}
\put(295,173){\vector(1,0){63}}
\put(23,173){\vector(1,0){28}}
\put(147,173){\vector(1,0){123}}
\put(23,123){\vector(1,0){28}}
\put(70,123){\vector(1,0){44}}
\put(164,123){\vector(1,0){106}}
\put(298,123){\vector(1,0){60}}
\put(160,73){\line(1,0){102}}
\put(160,71){\line(1,0){102}}
\end{picture}

\bigskip

On en d\'eduit une suite exacte
$$0\lra\E\oplus\ker(\pi)\lra\F'\oplus\F_0\lra{\cal U}\lra 0.$$
Ceci d\'emontre 1-.

Supposons maintenant que les hypoth\`eses de 2- soient v\'erifi\'ees.
Soient
$${\cal B} = (\Gamma\ot\Hom(\Gamma,\F))\oplus\F'.$$
On consid\`ere le morphisme injectif
$$B : \E'\lra{\cal B}$$
dont la premi\`ere composante est la compos\'ee
$$\E'\hfl{ev^*}{}\Gamma\ot\Hom(\E',\Gamma)^*\hfl{\lambda}{}
\Gamma\ot\Hom(\Gamma,\F)$$
et dont la seconde provient de \m{\Phi}. \Oad

\bigskip

\begin{picture}(360,230)
\put(135,220){$0$}
\put(280,220){$0$}
\put(10,170){$0$}
\put(55,170){$\E'$}
\put(133,170){$\cal A$}
\put(266,170){$\F_0\oplus\F'$}
\put(365,170){$0$}
\put(10,120){$0$}
\put(55,120){$\E'$}
\put(133,120){$\cal B$}
\put(84,125){$B$}
\put(265,120){$\coker(B)$}
\put(365,120){$0$}
\put(110,70){$\Gamma\ot\coker(\lambda)$}
\put(255,70){$\Gamma\ot\coker(\lambda)$}
\put(135,20){$0$}
\put(280,20){$0$}

\put(137,214){\vector(0,-1){30}}
\put(137,164){\vector(0,-1){30}}
\put(137,114){\vector(0,-1){30}}
\put(137,64){\vector(0,-1){25}}
\put(282,214){\vector(0,-1){30}}
\put(282,164){\vector(0,-1){30}}
\put(282,114){\vector(0,-1){30}}
\put(282,64){\vector(0,-1){30}}
\put(56,163){\line(0,-1){29}}
\put(58,163){\line(0,-1){29}}

\put(70,173){\vector(1,0){55}}
\put(310,173){\vector(1,0){51}}
\put(23,173){\vector(1,0){28}}
\put(147,173){\vector(1,0){114}}
\put(23,123){\vector(1,0){28}}
\put(70,123){\vector(1,0){55}}
\put(147,123){\vector(1,0){112}}
\put(318,123){\vector(1,0){40}}
\put(178,73){\line(1,0){72}}
\put(178,71){\line(1,0){72}}
\end{picture}

\bigskip

Puisque \ \m{\Ext^1(\Gamma,\F_0) = \Ext^1(\Gamma,\F') =
\lbrace 0\rbrace}, on a un isomorphisme
$$\coker(B) \ \simeq \ (\Gamma\ot\coker(\lambda))\oplus\F_0\oplus\F'.$$
Le carr\'e commutatif
$$\diagram{
\E'         & \hfl{B}{}    & {\cal B}                 \cr
\vfl{}{}    &              & \vfl{}{ev\oplus I_{\F'}} \cr
\E\oplus\E' & \hfl{\Phi}{} & \F\oplus\F' \cr
}$$
induit un morphisme surjectif
$$\rho : \coker(B)\lra\coker(\Phi) = {\cal U},$$
et une suite exacte
$$0\lra\E_0\lra\ker(\rho)\lra\E\lra 0.$$
Comme \ \m{\Ext^1(\E,\E_0)=\lbrace 0\rbrace}, on a un isomorphisme
$$\ker(\rho) \ \simeq \ \E\oplus\E_0.$$
On a donc une suite exacte
$$0\lra\E\oplus\E_0\lra
(\Gamma\ot\coker(\lambda))\oplus\F_0\oplus\F'\lra 0.$$
Ceci d\'emontre 2-. \carre

\bigskip

Les cas particuliers \ \m{\E_0 = 0} \ ou \m{\F_0 = 0} suivants seront
utilis\'es par la suite. Soit \m{M} un espace vectoriel de dimension finie.

\begin{xcoro}
1 - On suppose que \ \m{\Ext^1(\E,\Gamma) = \lbrace 0\rbrace}. Soient
$$\Phi : \E\oplus\E'\lra(\Gamma\ot M)\oplus\F'$$
un morphisme injectif induisant une surjection
$$\lambda : \Hom(\E',\Gamma)^*\lra M,$$
et \ \m{{\cal U}=\coker(\Phi)}. Alors il existe une
suite exacte
$$0\lra\E\oplus(\Gamma\ot\ker(\lambda))\lra\F_0\oplus\F'\lra
{\cal U}\lra 0.$$

\medskip

\noindent 2- On suppose que \ \m{\Ext^1(\Gamma,\F') = \lbrace 0\rbrace}.
Soient$$\Phi : \E\oplus(\Gamma\ot M)\lra\F\oplus\F'$$
un morphisme injectif induisant une injection
$$\lambda : M\lra\Hom(\Gamma,\F),$$
et \ \m{{\cal U}=\coker(\Phi)}.
Alors il existe une
suite exacte
$$0\lra\E\oplus\E_0\lra(\Gamma\ot\coker(\lambda))\oplus\F'\lra
{\cal U}\lra 0.$$
\end{xcoro}

\subsection{Applications}

Soient \m{X} une vari\'et\'e projective, \m{r,s} des entiers positifs, et
\m{\E_1,\ldots,\E_r,},\m{\F_1,\ldots,\F_s} des
faisceaux coh\'erents simples sur \m{X} tels que
$$\Hom(\E_i,\E_{i'}) = 0 \ \ {\rm si \ } i > i' \ , \
\Hom(\F_j,\F_{j'}) = 0 \ \ {\rm si \ } j > j', $$
$$\Hom(\F_j,\E_i) = \lbrace 0 \rbrace \ \
{\rm pour \ tous \ } i,j .$$
On suppose que pour \m{1\leq i\leq r} le morphisme canonique
$$\E_i\lra \Hom(\E_i,\F_1)^*\ot \F_1$$
est injectif. Soit \m{\G_i} son conoyau. Du corollaire 5.2, 1-, on
d\'eduit la

\begin{xprop}
1 - Soient
$$\Phi  : \som_{1\leq i\leq r}(\E_i\ot M_i) \lra
\som_{1\leq l\leq s}(\F_l\ot N_l)$$
un morphisme injectif, \m{\cal U} son conoyau
et \m{p} un entier tel que \m{0\leq p\leq r-1}. On suppose que
$$\Ext^1(\G_j,\F_l) = \Ext^1(\E_i,\F_1) =
\lbrace 0\rbrace$$
pour \m{p+1\leq j\leq r}, \m{1\leq i\leq p} et \m{2\leq l\leq s}.
Soit
$$f_p : \som_{p+1\leq j\leq r}(\Hom(\E_j,\F_1)^*\ot M_j)\lra N_1$$
l'application lin\'eaire d\'eduite de $\Phi$. On suppose que $f_p$ est
surjective. Alors il existe une suite exacte
$$0\lra \biggl(\bigoplus_{1\leq i\leq p}(\E_i\ot M_i)\biggr)\oplus
(\F_1\ot\ker(f_p))\lra
\biggr(\bigoplus_{p<j\leq r}(\G_j\ot M_j)\biggl)\oplus
\biggr(\bigoplus_{2\leq l\leq s}(\F_l\ot N_l)\biggl)\lra{\cal U}\lra 0.$$
\end{xprop}

\bigskip

Du corollaire 5.2, 2-, on d\'eduit la

\begin{xprop}
Soient \m{P_1} un espace vectoriel de dimension finie,
$$\Psi : \biggl(\bigoplus_{1\leq i\leq p}(\E_i\ot M_i)\biggr)\oplus
(\F_1\ot P_1)\lra
\biggr(\bigoplus_{p<j\leq r}(\G_j\ot M_j)\biggl)\oplus
\biggr(\bigoplus_{2\leq l\leq s}(\F_l\ot N_l)\biggl)$$
un morphisme injectif et \m{\cal U} son conoyau. On suppose que
$$\Ext^1(\F_1,\E_i) = \Ext^1(\E_i,\E_j) = \Ext^1(\F_1,\F_l) =
\lbrace 0\rbrace$$
pour \m{1\leq i\leq p}, \m{p+1\leq j\leq r}, \m{2\leq l\leq s}. Soit
$$g : P_1\lra \bigoplus_{p+1\leq j\leq r}\biggl(
\Hom(\F_1,\G_j)\ot M_j\biggr)$$
l'application lin\'eaire d\'eduite de \m{\Psi}. On suppose \m{g} injective.
Alors il existe une suite exacte
$$0\lra\som_{1\leq i\leq r}(\E_i\ot M_i) \lra (\F_1\ot\coker(g))\oplus
\biggl(\som_{2\leq l\leq s}(\F_l\ot N_l)\biggr)\lra{\cal U}\lra 0.$$
\end{xprop}

\vfil
\eject

\section{Mutations abstraites}

\subsection{Espaces abstraits de morphismes}

\subsubsection{D\'efinition g\'en\'erale}

Soient \m{\X{1}}, \m{\X{2}}, \m{\X{3}}, \m{\X{4}}, \m{M} , \m{H_L}, \m{H_R}
des espaces vectoriels sur un corps commutatif \m{k}, de dimension finie,
avec
$$\dim(M) < \dim(\X{2}).$$
On pose
$$W = (\X{1}\ot M)\oplus(\X{2}\ot M)\oplus \X{3}\oplus\X{4}.$$

Soient \m{\GG{0}}, \m{\GG{1}}, \m{\GG{2}} des groupes. On suppose que :
\begin{itemize}
\item[] $\GG{0}$ op\`ere lin\'eairement \`a gauche sur \m{\X{3}}, \m{\X{4}},
\m{H_R}.
\item[] $\GG{1}$ op\`ere lin\'eairement \`a droite sur \m{\X{1}}, \m{\X{3}},
\m{H_L}.
\item[] $\GG{2}$ op\`ere lin\'eairement \`a droite sur \m{\X{2}}, \m{\X{4}},
et \`a gauche sur \m{H_L}.
\end{itemize}
On suppose que ces actions sont {\em compatibles}, c'est-\`a-dire que
si deux de ces groupes \m{G_\alpha}, \m{G_\beta}, op\`erent sur un
m\^eme espace vectoriel \m{Z}, \`a gauche et \`a droite respectivement,
on a, pour tous \m{g_\alpha\in G_\alpha},\m{g_\beta\in G_\beta} et
\m{z\in Z},
$$g_\alpha(zg_\beta) = (g_\alpha z)g_\beta.$$
On suppose aussi que le groupe \m{\lbrace1,-1\rbrace} est contenu dans
\m{\GG{1}}, \m{\GG{2}} et \m{\GG{3}}, et agit comme on le pense sur les
espaces vectoriels sur lesquels ces groupes agissent (c'est-\`a-dire
que \m{-1} agit par multiplication par \m{-1}).

Soient
$$\g{3} : H_R\ot \X{1}\lra \X{3},$$
$$\g{4} : H_R\ot \X{2}\lra \X{4},$$
$$\g{1} : \X{2}\ot H_L\lra \X{1},$$
$$\g{2} : \X{4}\ot H_L\lra \X{3}$$
des applications lin\'eaires. On suppose que le diagramme suivant \m{(D)}
est commutatif :
$$\diagram{
H_R\ot \X{2}\ot H_L & \hfl{I_{H_R}\ot \g{1}}{} & H_R\ot \X{1}\cr
\vfl{\g{4}\ot I_{H_L}}{} & & \vfl{}{\g{3}} \cr
\X{4}\ot H_L & \hfl{\g{2}}{} & \X{3}
}$$
On suppose aussi que ces applications lin\'eaires sont compatibles avec
l'action des groupes. Par exemple \m{\GG{1}} op\`ere \`a droite sur \m{\X{1}}
et \m{H_L}, donc pour tous \m{\gg{1}\in \GG{1}}, \m{h_L\in H_L} et \m{y_1\in
\X{2}}
on a
$$\g{1}(y_1\ot (h_L\gg{1})) = \g{1}(y_1\ot h_L).\gg{1}.$$
De m\^eme, \m{\GG{2}} op\`ere \`a droite sur \m{\X{2}} et \`a gauche sur
\m{H_L}, donc pour tous \m{\gg{2}\in \GG{2}}, \m{h_L\in H_L} et \m{y_1\in
\X{2}}
on a
$$\g{1}(y_1\gg{2}\ot h_L) = \g{1}(y_1\ot \gg{2}h_L).$$
On suppose aussi que \m{\g{4}} est surjective, et que l'application
lin\'eaire
$$\ov{\g{1}} : H_L\lra \X{2}^*\ot \X{1}$$
d\'eduite de \m{\g{1}} est injective.

\begin{defin}
On appelle {\em espace abstrait de morphismes} et on note \m{\Theta}
la donn\'ee de \m{\X{1}}, \m{\X{3}}, \m{\X{2}}, \m{\X{4}}, \m{H_L}, \m{H_R},
\m{\GG{0}}, \m{\GG{1}}, \m{\GG{2}}, \m{\g{1}}, \m{\g{2}}, \m{\g{3}} et
\m{\g{4}}.
L'espace vectoriel
$$W = (\X{1}\ot M)\oplus(\X{2}\ot M)\oplus \X{3}\oplus \X{4}$$
est {\em l'espace total} de \m{\Theta}.
\end{defin}

\subsubsection{Dictionnaire}

Si on \'etudie les morphismes
$$\E\oplus\E'\lra(\Gamma\ot M)\oplus\F',$$
l'espace abstrait de morphismes associ\'e est d\'efini par
$$\X{1} = \Hom(\E,\Gamma), \ \X{2} = \Hom(\E',\Gamma),$$
$$\X{3} = \Hom(\E,\F'), \ \X{4} = \Hom(\E',\F'),$$
$$H_L = \Hom(\E,\E'), \ H_R = \Hom(\Gamma,\F'),$$
$$\GG{0} = Aut(\F'), \ \GG{1} = Aut(\E), \ \GG{2} = Aut(\E'),$$
les applications \m{\g{1}}, \m{\g{2}}, \m{\g{3}}, \m{\g{4}}
\'etant les compositions des morphismes. On a dans ce cas
$$W = \Hom(\E\oplus\E',(\Gamma\ot M)\oplus\F').$$

\subsection{Groupes associ\'es}

On va construire deux nouveaux groupes associ\'es \`a \m{\Theta} :
\m{G_L} et \m{G_R}. Le groupe \m{G_L} est constitu\'e des matrices
$$\pmatrix{
\gg{1} & 0\cr h_L & \gg{2}}$$
avec \m{\gg{1}\in \GG{1}}, \m{\gg{2}\in \GG{2}}, \m{h_L\in H_L}.
La loi de groupe de \m{G_L} est
$$\pmatrix{\gg{1} & 0\cr h_L & \gg{2}}.\pmatrix{\gg{1}' & 0\cr h_L' &
 \gg{2}'} =
\pmatrix{\gg{1}\gg{1}' & 0\cr h_L\gg{1}'+ \gg{2}h_L' & \gg{2}\gg{2}'}.$$
Le groupe \m{G_R} est constitu\'e des matrices
$$\pmatrix{
g_M & 0\cr \lambda & g_0}$$
avec \m{g_M\in GL(M)}, \m{g_0\in \GG{0}}, \m{\lambda\in M^*\ot H_R}.
La loi de groupe de \m{G_R} est
$$\pmatrix{g_M & 0\cr \lambda & g_0}.\pmatrix{g_M' & 0\cr \lambda' & g_0'} =
\pmatrix{g_Mg_M' & 0\cr \lambda g_M'+ g_0\lambda' & g_0g_0'}$$
(\m{GL(M)} agit de mani\`ere \'evidente \`a droite sur le premier facteur
de \m{M^*\ot H_R}, et \m{\GG{0}} \`a gauche sur le deuxi\`eme facteur).

Dans le cas du \para 6.1, on a
$$G_L=Aut(\E\oplus\E'), \ G_R=Aut((\Gamma\ot M)\oplus\F').$$

\subsection{Actions des groupes associ\'es sur l'espace de morphismes}

Le groupe \m{G_L} op\`ere \`a droite sur \m{W} : si \m{\phi_1\in \X{1}\ot M},
\m{\psix\in \X{2}\ot M}, \m{x_2\in \X{3}}, \m{y_2\in \X{4}}, \m{\gg{1}\in
\GG{1}},
\m{\gg{2}\in \GG{2}} et \m{h_L\in H_L} on a
$$\pmatrix{\phi_1 &\psix\cr x_2 & y_2}\pmatrix{\gg{1} & 0\cr h_L & \gg{2}} =
\pmatrix{\phi_1\gg{1}+(\g{1}\ot I_M)(\psix\ot h_L) & \psix \gg{2} \cr
x_2\gg{1} + \g{2}(y_2\ot h_L) & y_2\gg{2}}.$$

Le groupe \m{G_R} op\`ere \`a gauche sur \m{W} : si \m{\phi_1\in \X{1}\ot M},
\m{\psix\in \X{2}\ot M}, \m{x_2\in \X{3}}, \m{y_2\in \X{4}}, \m{g_0\in \GG{0}},
\m{g_M\in GL(M)} et \m{\lambda\in M^*\ot H_R} on a
$$\pmatrix{g_M & 0\cr \lambda & g_0}\pmatrix{\phi_1 &\psix\cr x_2 & y_2} =
\pmatrix{(I_{\X{1}}\ot g_M)(\phi_1) & (I_{\X{2}}\ot g_M)(\psix)\cr
g_0x_2+\g{3}(\pline{\lambda,\phi_1}) &
g_0y_2+\g{4}(\pline{\lambda,\psix})}.$$

Les actions de ces groupes sont compatibles, c'est-\`a-dire que si
\m{g_L\in G_L}, \m{g_R\in G_R} et \m{w\in W}, on a
$$g_R(wg_L) = (g_Rw)g_L.$$
C'est pourquoi on parlera abusivement du groupe \ \m{G_L\times G_R}
\ ou d'un de ses sous-groupes et de son action sur \m{W} (au lieu
d'utiliser par exemple le groupe \m{G_L^{op}\times G_R}).

On note \m{H} le \og sous-groupe\fg \ de \ \m{G_L\times G_R} constitu\'e
des paires
$$\pmatrix{\pmatrix{1 & 0 \cr h_L & 1},
\pmatrix{1 & 0 \cr \lambda & 1}},$$
(o\`u \ \m{h_l\in H_L}, \m{\lambda\in M^*\ot H_R}).

\subsection{Mutation d'un espace abstrait de morphismes}

Soit \m{M'} un \m{k}-espace vectoriel tel que
$$\dim(M') = \dim(\X{2})-\dim(M).$$
On va d\'efinir un nouvel espace abstrait de morphismes \m{D(\Theta)}
associ\'e \`a \m{\Theta}. Posons
$$\X{1}' = H_R, \ \ \X{2}' = \X{2}^*, \ \ \X{3}' = \X{3}, \ \
\X{4}' = \coker(\ov{\g{1}}) = (\X{2}^*\ot\X{1})/H_L,$$
$$H'_R = \X{1}, \ \ H'_L = \ker(\g{4})\subset H_R\ot\X{2}.$$
Soient
$$\g{1}' : \X{2}'\ot H'_L\lra \X{1}'$$
la restriction de la contraction
$$\X{2}^*\ot\X{2}\ot H_R\lra H_R,$$
$$\g{3}' = \g{3} : H'_R\ot\X{1}'\lra\X{3}',$$
et
$$\g{4}' : H'_R\ot\X{2}'\lra\X{4}'$$
la projection
$$\X{2}^*\ot\X{1}\lra(\X{2}^*\ot\X{1})/H_L.$$
La d\'efinition de \m{\g{2}'} est un peu plus compliqu\'ee. On a un
diagramme commutatif, o\`u la ligne du haut et la colonne de gauche sont
commutatives :

\bigskip
\bigskip

\begin{picture}(360,230)
\put(90,220){$0$}
\put(10,170){$0$}
\put(55,170){$\ker(\g{4})\ot H_L$}
\put(168,170){$\ker(\g{4})\ot\X{2}^*\ot \X{1}$}
\put(301,170){$\ker(\g{4})\ot\coker(\ov{\g{1}})$}
\put(430,170){$0$}
\put(55,120){$H_R\ot\X{2}\ot H_L$}
\put(176,120){$H_R\ot\X{2}\ot\X{2}^*\ot \X{1}$}
\put(75,70){$\X{4}\ot H_L$}
\put(188,70){$\X{3}$}
\put(90,20){$0$}

\put(92,214){\vector(0,-1){30}}
\put(92,164){\vector(0,-1){30}}
\put(92,114){\vector(0,-1){30}}
\put(92,64){\vector(0,-1){25}}
\put(192,164){\vector(0,-1){30}}
\put(194,94){$\phi$}
\put(192,114){\vector(0,-1){30}}

\put(133,173){\vector(1,0){30}}
\put(410,173){\vector(1,0){15}}
\put(23,173){\vector(1,0){28}}
\put(273,173){\vector(1,0){22}}
\put(138,129){$I\ot\ov{\g{1}}$}
\put(133,123){\vector(1,0){40}}
\put(136,78){$\g{2}$}
\put(128,72){\vector(1,0){55}}
\end{picture}

\bigskip
Le morphisme \m{\phi} est la contraction de \m{\X{2}\ot\X{2}^*}, suivie
de \m{\g{3}}. La commutativit\'e du carr\'e du bas d\'ecoule de celle du
carr\'e \m{(D)} du \para 6.1.1. Il en d\'ecoule que \m{\phi} s'annule
sur \ \m{\ker(\g{4})\ot H_L}, et induit donc une application lin\'eaire
$$\g{2}' : \X{4}'\ot H_L' = \coker(\ov{\g{1}})\ot\ker(\g{4})
\lra\X{3}=\X{3}'.$$
Il est ais\'e de voir que l'analogue du carr\'e \m{(D)} du \para 6.1.1
est commutatif. Il est clair que \m{\g{1}'} induit une injection
$$\ov{\g{1}'} : H'_L\lra\X{2}'^*\ot\X{1}',$$
(c'est l'inclusion \ \m{\ker(\g{4})\subset H_R\ot\X{2}}), et que \m{\g{4}'}
est surjective. On pose
$$\GG{0}' = \GG{1}^{op}, \ \GG{1}' = \GG{0}^{op}, \ \GG{2}' = \GG{2}^{op}.$$
Les actions de ces groupes se d\'eduisent imm\'ediatement de celles des
groupes \m{\GG{0}},\m{\GG{1}} et \m{\GG{2}}. Par exemple, \m{\GG{1}} agit \`a
droite
sur \m{H_L} et \m{\X{1}}, et cette action est compatible avec
$$\g{1} : \X{2}\ot H_L\lra\X{1}.$$
On obtient donc une action \`a droite de \m{\GG{1}} sur
\m{(\X{2}^*\ot\X{1})/H_L}, c'est-\`a-dire une action \`a gauche de
\m{\GG{0}'} sur \m{\X{4}'}.

\begin{defin}
On note \m{D(\Theta)} l'espace abstrait de morphismes d\'efini par
\m{\X{1}'}, \m{\X{2}'}, \m{\X{3}'}, \m{\X{4}'}, \m{H'_L}, \m{H'_R},
\m{\GG{0}'}, \m{\GG{1}'}, \m{\GG{2}'}, \m{\g{1}'}, \m{\g{2}'}, \m{\g{3}'} et
\m{\g{4}'}. On l'appelle la {\em mutation de \m{\Theta}}.
\end{defin}

\begin{xprop}
On a \ \ \m{D(D(\Theta)) = \Theta}.
\end{xprop}

Imm\'ediat. \carre

On d\'efinit comme pour \m{\Theta} les \og groupes\fg \ \m{G'_L\times G'_R}
et \m{H'} correspondant \`a \m{D(\Theta)}.

\subsection{Mutation des morphismes}

On note \m{W'} l'espace total de \m{D(\Theta)}, c'est-\`a-dire
$$W' = (\X{1}'\ot M')\oplus (\X{2}'\ot M') \oplus\X{3}'\oplus\X{4}'.$$
On note \m{W^0} l'ouvert de \m{W} constitu\'e des
$$\pmatrix{\phi_1 & \psix\cr \x{3} & \x{4}}$$
tels que l'application lin\'eaire
$$\ov{\psix} : \X{2}^*\lra M$$
d\'eduite de \m{\psix} soit surjective. On d\'efinit de m\^eme l'ouvert
\m{W'^0} de \m{W'}.

Rappelons que la projection
$$\X{2}^*\ot\X{1}\lra(\X{2}^*\ot\X{1})/H_L$$
n'est autre que \m{\g{4}'}. De m\^eme, la projection
$$\X{2}'^*\ot\X{1}'\lra({\X{2}'}^*\ot\X{1}')/H'_L$$
n'est autre que \m{\g{4}}. Si \ \m{\psix\in\X{2}\ot M}, on notera
\m{q(\psix)} l'application lin\'eaire
$$\ov{\psix}\ot I_{\X{1}} : \X{2}^*\ot\X{1}\lra M\ot\X{1}.$$
On d\'efinit de m\^eme, pour tout \ \m{\psix'\in\X{2}'\ot M'}
l'application lin\'eaire
$$q'(\psix') : \X{2}\ot H_R = {\X{2}'}^*\ot\X{1}'\lra M'\ot\X{1}'.$$

Soit
$$w = \pmatrix{\phi_1 & \psix\cr \x{3} & \x{4}}\in W^0.$$
On va en d\'eduire un \'el\'ement de \m{W'^0} (pas de mani\`ere unique).
On choisit d'abord un isomorphisme
$$\ker(\ov{\psix})^*\simeq M'.$$
On note \m{\psix'} l'\'el\'ement de \ \m{\X{2}'\ot M'} \ provenant de
l'application lin\'eaire
$$\ov{\psix'} : {\X{2}'}^* = \X{2}\lra\ker(\ov{\psix})^*=M',$$
qui est la transpos\'ee de l'inclusion de \m{\ker(\ov{\psix})} dans
\m{\X{2}^*}.

Soit
$$u\in\g{4}^{-1}(-\x{4})\subset H_R\ot\X{2}.$$
Notons que \m{u} est d\'efini \`a un \'el\'ement pr\`es de
\ \m{\ker(\g{4})=H'_L}.
Soit
$$\phi'_1 = q'(\psix')(u) \in \X{1}'\ot M'.$$
On peut aussi \'ecrire
$$\phi'_1 = \pline{\psix',u}. $$

Soient
$$\alpha\in q(\psix)^{-1}(\phi_1)\subset \X{2}^*\ot\X{1},$$
et
$$\x{4}' = \g{4}'(\alpha)\in\X{4}'.$$
Notons que \m{\alpha} est d\'efini \`a un \'el\'ement pr\`es de
$$\ker(\ov{\psix})\ot\X{1} = H'_R\ot {M'}^*.$$

Soient enfin
$$\x{3}' = \x{3}+\g{3}(\pline{\alpha,u}).$$
et
$$z(w,u,\alpha) = \pmatrix{\phi'_1 & \psix'\cr \x{3}' & \x{4}'}
\in W'^0.$$
On emploiera aussi la notation
$$z(w) = z(w,u,\alpha)$$
bien que cet \'el\'ement de \m{W'^0} ne d\'epende pas uniquement de \m{w}.

\begin{xprop}
Soit \m{w\in W^0}. Les \'el\'ements \m{z(w,u,\alpha)}, pour tous les choix
possibles de \m{u} et \m{\alpha}, constituent une \m{H'}-orbite de \m{W'^0}.
\end{xprop}

\dem On v\'erifie ais\'ement que si on remplace
\m{u} par \m{u+h'_L} et \m{\alpha} par \m{\alpha+\psi} (avec
\m{h'_L\in H'_L} et \m{\psi\in M^*\ot H_R}), l'\'el\'ement obtenu de
\m{W'} est
$$\pmatrix{1 & 0 \cr \psi & 1}\pmatrix{\phi_1' & \psix' \cr \x{3}' &
\x{4}'}\pmatrix{1 & 0 \cr h'_L & 1}.$$
\carre

\begin{xprop}
Pour tout \m{w\in W^0}, on a
$$z(z(w)) \ \in \ (G_L\times G_R)w.$$
\end{xprop}

\dem On part de
$$w = \pmatrix{\phi_1 & \psix \cr \x{3} & \x{4}} \in W^0,$$
et on prend comme pr\'ec\'edemment \ \m{u\in\g{4}^{-1}(-\x{4})},
\m{\alpha\in q(\psix)^{-1}(\phi_1)} \ pour d\'efinir
$$z(w) = \pmatrix{\phi_1 & \psix' \cr \x{3}' & \x{4}'} \in {W'}^0.$$
On cherche maintenant \ \m{u'\in\g{4}'^{-1}(-\x{4}')} \ et \
\m{\alpha\in q(\psix')^{-1}(\phi'_1)} \ pour d\'efinir \m{z(z(w))}.
On a
$$\g{4}'(-\alpha) = - \x{4}',$$
donc on peut prendre
$$u' = - \alpha.$$
D'autre part, on a
$$q'(\psix') = \phi'_1,$$
et on peut prendre
$$\alpha' = u.$$
Soit
$$z(z(w)) = \pmatrix{\phi_1 & \psix'' \cr \x{3}'' & \x{4}''} \in {W}^0$$
l'\'el\'ement de \m{W^0} d\'efini par \m{u'} et \m{\alpha}. On a
\'evidemment \ \m{\psix'' = \psix}, et
$$\phi_1'' = q(\psix)(u') = -q(\psix)(\alpha) = -\phi_1,$$
$$\x{4}'' = \g{4}(\alpha') = \g{4}(u) = -\x{4},$$
$$\x{3}'' = \x{3}' + \g{3}(\pline{\alpha',u'}) =
\x{3} + \g{3}(\pline{\alpha,u}) - \g{3}(\pline{\alpha,u}) = \x{3}.$$
Donc
\begin{eqnarray*}
z(z(w)) & = & \pmatrix{-\phi_1 & \psix \cr \x{3} & -\x{4}} \cr
        & = & \pmatrix{-1 & 0 \cr 0 & 1}
\pmatrix{\phi_1 & \psix \cr \x{3} & \x{4}}\pmatrix{1 & 0 \cr 0 & -1}\cr
\end{eqnarray*}
\carre

\begin{xprop}
Soient \m{w_1,w_2\in W^0} des points qui sont dans la m\^eme \
\m{(G_L\times G_R)}-orbite. Alors \m{z(w_1)} et \m{z(w_2)} sont dans la
m\^eme \ \m{(G'_L\times G'_R)}-orbite.
\end{xprop}

\dem On v\'erifie ais\'ement que c'est vrai si
$$w_2 = \pmatrix{g_M & 0 \cr 0 & g_0}w_1\pmatrix{g_1 & 0 \cr 0 & g_2}.$$
Il reste \`a traiter les cas
$$w_2 = \pmatrix{1 & 0 \cr \psi & 1}w_1,$$
ou
$$w_2 = w_1\pmatrix{1 & 0 \cr h_l & 1},$$
avec \m{h_l\in H_L}, \m{\psi\in M^*\ot H_R}. On ne traitera que le premier
cas, le second \'etant analogue. Posons
$$w_1 = \pmatrix{\phi_1 & \psix \cr \x{3} & \x{4}}.$$
Alors on a
$$w_2 = \pmatrix{\phi_1 + (\g{1}\ot I_M)(\psix\ot h_L) & \psix \cr
\x{3} + \g{2}(\x{4}\ot h_L) & \x{4}}.$$
On suppose que
$$z(w_1) = \pmatrix{\phi'_1 & \psix' \cr \x{3}' & \x{4}'}$$
est d\'efini par \ \m{u_1\in H_R\ot X_2} \ et \
\m{\alpha_1\in\X{2}^*\ot\X{1}}. On va chercher des \'el\'ements \m{u_2},
\m{\alpha_2} convenables pour d\'efinir \m{z(w_2)}. On doit avoir
\ \m{u_2\in\g{4}^{-1}(-\x{4})}, donc on peut prendre
$$u_2 = u_1.$$
On doit avoir
$$q(\psix)(\alpha_2) = \phi_1 + (\g{1}\ot I_M)(\psix\ot h_L).$$
Posons
$$\alpha_2 = \alpha_1+\alpha_0.$$
On doit donc avoir
$$q(\psix)(\alpha_0) = (\g{1}\ot I_M)(\psix\ot h_L).$$
Pour cela il suffit de prendre
$$\alpha_0 = h_L$$
(vu comme \'el\'ement de \ \m{\X{2}^*\ot\X{1}}, \`a l'aide de
\m{\ov{\g{1}}}). On a alors
$$z(w_2) = \pmatrix{\phi_1'' & \psix' \cr \x{3}'' & \x{4}''},$$
avec
$$\phi_1'' = q'(\psix')(u_2) = q'(\psix')(u_1) = \phi_1',$$
$$\x{4}'' = \g{4}(\alpha_2) = \g{4}(\alpha_1+\alpha_0)
= \g{4}(\alpha_1) = \x{4}',$$
\begin{eqnarray*}
\x{3}'' & = & \x{3}+\g{3}(\pline{\alpha_2,u_2}) \cr
        & = & \x{3}'-\g{2}(h_L\ot\g{4}(u_2))+
\g{3}(\pline{\alpha_0,u_2})\cr
\end{eqnarray*}
Mais on a
$$\g{2}(h_L\ot\g{4}(u_2))=\g{3}(\pline{\alpha_2,u_2}),$$
si on se souvient que \ \m{\alpha_0=\ov{\g{1}}(h_L)}, \`a cause du
diagramme commutatif \m{(D)} du \para 6.1.1. On a donc
\ \m{\x{3}''=\x{3}'} \ et finalement
$$z(w_1)=z(w_2).$$
\carre

\subsection{Th\'eor\`emes d'isomorphisme}

Le th\'eor\`eme suivant d\'ecoule imm\'ediatement des r\'esultats du
\para 6.5 :

\begin{xtheo}
L'application associant \`a l'orbite d'un point \m{w} de \m{W^0} l'orbite
de \m{z(w)} d\'efinit une bijection
$$D_\Theta : W^0/(G_L\times G_R) \ \simeq \ {W'}^0/(G_L'\times G_R').$$
\end{xtheo}

\bigskip

On suppose maintenant que les groupes \m{\GG{1}}, \m{\GG{2}}, \m{\GG{3}},
sont alg\'ebriques sur \m{k} et que leurs actions sont alg\'ebriques. Il
est alors clair par construction que pour tout \m{w\in W^0}, il
existe un voisinage de Zariski \m{U} de \m{w} dans \m{W^0} et un
voisinage de Zariski \m{U'} de \m{z(w)} dans \m{{W'}^0} tels que
\m{D_\Theta} se rel\`eve en un morphisme \ \m{U\lra U'} et que
\m{D_{D(\Theta)}} se rel\`eve en un morphisme \ \m{U'\lra U}. On en
d\'eduit ais\'ement le r\'esultat suivant :

\begin{xtheo}
Soit \m{U} un ouvert \m{(G_L\times G_R)}-invariant de \m{W^0} tel qu'il
existe un bon quotient \m{U//(G_L\times G_R^{op})}. Soit \m{U'}
l'ensemble des points de \m{{W'}^0} au dessus de

\noindent
\m{D_\Theta(U/((G_L\times G_R))}. Alors \m{U'} est un ouvert
\m{(G'_L\times G'_R)}-invariant de \m{{W'}^0}, et il existe un bon
quotient \m{U'//(G'_L\times {G'_R}^{op})}, qui est isomorphe \`a
\m{U//(G_L\times G_R^{op})}.
\end{xtheo}

\bigskip
\bigskip

\section{Mutations de morphismes de type \m{(r,s)}}

\subsection{Espaces de morphismes abstraits associ\'es}

On applique les r\'esultats du \para 6 aux cas d\'ecrits au \para 3.
Commen\c cons par d\'ecrire la situation en termes de morphismes de
faisceaux. On s'int\'eresse aux morphismes
$$\som_{1\leq i\leq r}(\E_i\ot M_i)\lra\som_{1\leq l\leq s}(\F_l\ot N_l).$$
Soit $p$ un entier tel que \m{0\leq p\leq r-1}. On pose
$$\E=\som_{1\leq i\leq p}(\E_i\ot M_i), \
\E'=\som_{p+1\leq j\leq r}(\E_j\ot M_j),$$
$$\Gamma=\F_1, \ M = N_1, \ \F'=\som_{2\leq l\leq s}(\F_l\ot N_l).,$$
de sorte que les morphismes pr\'ec\'edents peuvent s'\'ecrire sous la
forme
$$\E\oplus\E'\lra(\Gamma\ot M)\oplus\F',$$
comme dans le \para 5 et le \para 6.

On reprend maintenant le language du \para 3. On supposera que
$$n_1 < \sigg_{1\leq i\leq r}\dim(H_{1i})m_i.$$
On pose
$$X_1 = \som_{1\leq i\leq p}(H_{1i}\ot M_i^*), \
X_2 = \som_{p+1\leq j\leq r}(H_{1j}\ot M_j^*), $$
$$X_3 = \som_{1\leq i\leq p, 2\leq l\leq s}(H_{li}\ot M_i^*\ot N_l), \
X_4 = \som_{p+1\leq j\leq r, 2\leq l\leq s}(H_{lj}\ot M_j^*\ot N_l), $$
$$H_L=\som_{1\leq i\leq p,p+1\leq j\leq r}(A_{ji}\ot M_i^*\ot M_j), \
M=N_1, \ H_R=\som_{2\leq l\leq s}(B_{l1}\ot N_l).$$
On d\'efinit de m\^eme les groupes \m{\GG{0}}, \m{\GG{1}}, \m{\GG{2}}, et
les applications \m{\g{1}},\m{\g{2}},\m{\g{3}}, et \m{\g{4}}. On obtient
ainsi un espace abstrait de morphismes not\'e \m{\Theta_p}, d'espace
total de morphismes not\'e \m{W_p}. Il est clair que
$$W_p = W = \som_{1\leq i\leq r, 1\leq s\leq l}(H_{li}\ot M_i^*\ot N_l),$$
mais en g\'en\'eral
$$W_p^0\not = W_q^0$$
si \m{p\not = q}.

Soit
$$w=(\phi_{li})_{1\leq i\leq r,1\leq s\leq s}\in W.$$
Alors, par d\'efinition, on a \ \m{w\in W_p^0} \ si et seulement si
l'application lin\'eaire
$$\sigg_{p+1\leq j\leq r}\phi_j :
\som_{p+1\leq j\leq r}(H_{1j}^*\ot M_j)\lra N_1$$
est surjective. On a donc
$$W_{r-1}^0\subset W_{r-2}^0\subset\cdots\subset W_0^0\subset W.$$

\subsection{Description des mutations d'espaces de morphismes abstraits}

L'espace de morphismes abstraits \m{D(\Theta_p)} correspond \`a un
espace de morphismes de type \m{(r+s-p-1,p+1)}, dont nous allons d\'ecrire
le dual (cf. \para 3.6), qui est donc un espace de morphismes de type
\m{(p+1,r+s-p-1)}. Il est d\'efini par les espaces vectoriels
\m{M^{(p)}_1,\ldots,M^{(p)}_r},\m{N^{(p)}_1,\ldots,N^{(p)}_s},
\m{A^{(p)}_{ji},B^{(p)}_{ml},H^{(p)}_{li}}, et par les compositions
ad\'equates. On a
$$M^{(p)}_i=M_i \ \ {\rm si} \ 1\leq i\leq p,$$
$$\dim(M^{(p)}_{p+1}) = \biggl(\sigg_{p+1\leq j\leq r}m_j\dim(H_{1j})
\biggr)-n_1,$$
$$N^{(p)}_i=M_{p+i} \ \  {\rm si} \ 1\leq i\leq r-p,\ \ \
N^{(p)}_l = N_{l-r+p+1} \ \ {\rm si} \ r-p+1\leq l\leq r+s-p-1,$$
$$A^{(p)}_{ji}=A_{ji} \ \ {\rm si} \ 1\leq i\leq j\leq p,\ \ \
A^{(p)}_{p+1,i}=H_{1i} \ \ {\rm si} \ 1\leq i\leq p,$$
$$B^{(p)}_{lm} = A_{l+p,m+p} \ \ {\rm si} \ 1\leq m\leq l\leq r-p,$$
$$B^{(p)}_{lm} = B_{l-r+p+1,m-r+p+1} \ \ {\rm si} \ r-p+1\leq m\leq l
\leq r+s-p-1,$$
$$B^{(p)}_{lm} = \ker(B_{l-r+p+1,1}\ot H_{1,m+p}\lra H_{l-r+p+1,m+p})\ \ \ \
\ \ $$
$$\ \ \ \ \ \ \ \ {\rm si} \ r-p+1\leq l\leq r+s-p-1, \ 1\leq m\leq r - p,$$
$$H^{(p)}_{li} = (H_{1i}\ot H^*_{1,l+p})/A_{l+p,i} \ \ {\rm si} \
1\leq i\leq p, 1\leq l\leq r-p,$$
$$H^{(p)}_{li} = H_{l-r+p+1,i} \ \ {\rm si} \ 1\leq i\leq p,
r-p+1\leq l\leq r+s-p-1,$$
$$H^{(p)}_{l,p+1} = H^*_{1,l+p} \ \ {\rm si} \ 1\leq l\leq r-p,\ \ \
H^{(p)}_{l,p+1} = B_{l-r+p+1,1} \ \ {\rm si} \ r-p+1\leq l\leq r+s-p-1.$$
La description des compositions est laiss\'ee au lecteur. On note
$$W'(p) = \som_{1\leq i\leq p+1,1\leq l\leq r+s-p-1}
L(H^{(p)*}_{li}\ot M^{(p)}_i,N^{(p)}_l)$$
l'espace total de morphismes de \m{D(\Theta_p)}, \m{W'_0(p)} l'ouvert
correspondant, \m{G(p)} le groupe agissant sur \m{W'(p)}. Le th\'eor\`eme
6.5 dit qu'on a une bijection
$$W^0_p/G\ \simeq\ W'_0(p)/G(p).$$

\subsection{Description des mutations de morphismes}

Soit
$$w = (\phi_{li})_{1\leq i\leq r,1\leq l\leq s}\in W_p^0\subset
\som_{1\leq i\leq r,1\leq l\leq s}L(H_{li}^*\ot M_i,N_l),$$
et
$$\phi_1 = (\phi_{1i})_{1\leq i\leq p}, \ \ \phi_2 =
(\phi_{1j})_{p+1\leq j\leq r}, \ \ x_3 = (\phi_{li})_{1\leq i\leq p,
2\leq l\leq s}, \ \ x_4 = (\phi_{li})_{p+1\leq i\leq r,
2\leq l\leq s}.$$
On va d\'ecrire
$$z(w) = \pmatrix{\phi'_1 & \phi'_2\cr x'_3 & x'_4}.$$
On construit d'abord les \'el\'ements \m{u} et \m{\alpha} du \para 6.5.
On doit prendre pour \m{u} un \'el\'ement de
$$\som_{p+1\leq j\leq r,2\leq l\leq s}(H_{1j}\ot B_{l1}\ot M_j^*\ot N_l)$$
tel que le diagramme suivant soit commutatif :

\bigskip

\begin{picture}(360,100)
\put(100,90){$\som_{p+1\leq j\leq r,2\leq l\leq s}(H_{lj}^*\ot M_j)$}
\put(90,10){$\som_{p+1\leq j\leq r,2\leq l\leq s}
(H_{1j}^*\ot B_{l1}^*\ot M_j)$}
\put(320,90){$\som_{2\leq l\leq s}N_l$}

\put(225,92){\vector(1,0){85}}
\put(120,65){\vector(0,-1){40}}
\put(225,25){\vector(3,2){80}}
\put(245,95){$-x_4$}
\put(250,50){$u$}
\end{picture}

\bigskip

On prend pour \m{\alpha} un \'el\'ement de
$$\som_{1\leq i\leq p,p+1\leq j\leq r}(H_{1j}^*\ot H_{1i}\ot M_i^*\ot M_j)$$
tel que le diagramme suivant soit commutatif :

\bigskip

\begin{picture}(360,100)
\put(100,90){$\som_{1\leq i\leq p}(H_{li}^*\ot M_i)$}
\put(320,90){$N_1$}
\put(280,10){$\som_{p+1\leq j\leq r}(H_{lj}^*\ot M_j)$}

\put(200,92){\vector(1,0){110}}
\put(325,25){\vector(0,1){55}}
\put(170,75){\vector(3,-2){95}}
\put(220,45){$\alpha$}
\put(330,50){$\ov{\phi_2}$}
\end{picture}

\bigskip

Dans le premier diagramme, la fl\`eche verticale est induite par les
compositions \break \m{B_{l1}\ot H_{1j}\lra H_{lj}}, et dans le second, la
fl\`eche horizontale provient de \m{\phi_1}.

D\'eterminons maintenant \m{\phi'_1}, \m{\phi'_2}, \m{x'_3} et \m{x'_4}.
On prend \ \m{M'=\ker(\ov{\phi_2})}, qui est donc un quotient de
\ \m{\som_{p+1\leq j\leq r}(H_{lj}^*\ot M_j)}. Alors \m{x'_4} est l'image
de \m{\alpha} dans
$$\som_{1\leq i\leq p,p+1\leq j\leq r}\Biggl(\biggl((H_{1j}^*\ot H_{1i})
/A_{ji}\biggr)\ot M_i^*\ot M_j\Biggr)\ = \
\som_{1\leq i\leq p,1\leq l\leq r-p}(H^{(p)}_{li}\ot M^{(p)*}_i\ot
N^{(p)}_l),$$
\m{\phi'_1} est la restriction \`a \m{\ker(\ov{\phi_2})} de
l'aaplication lin\'eaire
$$\som_{p+1\leq j\leq r}(H_{1j}^*\ot M_j)\lra \som_{2\leq l\leq s}
(B_{l1}\ot N_l)$$
provenant de \m{u} et \m{\phi'_2} est l'inclusion
$$\ker(\ov{\phi_2})\subset \som_{p+1\leq j\leq r}(H_{lj}^*\ot M_j).$$
Pour obtenir \m{x'_3}, on fait la somme de \m{x_3} et de la compos\'ee
$$\som_{1\leq i\leq r}(H_{1i}^*\ot M_i)\hfl{\alpha}{}
\som_{p+1\leq j\leq r}(H_{1j}^*\ot M_j)\hfl{u}{}\som_{2\leq l\leq s}
(B_{l1}\ot N_l).$$

\subsection{Polarisations associ\'ees}

Soit \m{(\lambda_1,\ldots,\lambda_r,\mu_1,\ldots,\mu_s)} une polarisation
de l'action de \m{G} sur \m{W}. On va en d\'eduire une polarisation de
l'action de \m{G(p)} sur \m{W'(p)}.

Soient \ \m{M'_i\subset M_i}, \m{N'_l\subset N_l} ,\m{1\leq i\leq r},
\m{1\leq l\leq s} des sous-espaces vectoriels. On pose
$$m'_i=\dim(M'_i), \ \ n'_l=\dim(N'_l),$$
$${M^{(p)}}'_i=M'_i\ \ {\rm si} \ 1\leq i\leq p,$$
$${N^{(p)}}'_l=M'_{l-p} \ \ {\rm si} \ 1\leq l\leq r-p,\ \ \
{N^{(p)}}'_l=N'_{l-r+p+1} \ \ {\rm si} r-p+1\leq l\leq r+s-p-1,$$
$${m^{(p)}}'_i=\dim({M^{(p)}}'_i)\ \ {\rm si} \ 1\leq i\leq p, \ \ \
{m^{(p)}}'_{p+1} = \biggl(\sigg_{p+1\leq j\leq r}\dim(H_{1j})m'_j\biggr)-
n'_1,$$
$${n^{(p)}}'_l=\dim({N^{(p)}}'_l)\ \ {\rm si} \ 1\leq l\leq r+s-p-1.$$
On d\'efinit une suite
\m{(\alpha'_1,\ldots,\alpha'_{p+1},}\m{\beta'_1,\ldots,\beta'_{r+s-p-1})} de
nombres rationnels par les identit\'es
$$\sigg_{1\leq i\leq r}\lambda_im'_i-\sigg_{1\leq l\leq s}\mu_ln'_l =
\sigg_{1\leq i\leq p+1}\alpha'_i{m^{(p)}}'_i-\sigg_{1\leq l\leq r+s-p-1}
\beta'_l{n^{(p)}}'_l.$$
On a donc
$$\alpha'_i=\lambda_i\ \ {\rm si}\ 1\leq i\leq p,\ \ \
\alpha'_{p+1} = \mu_1,$$
$$\beta'_l=\mu_1\dim(H_{1,l+p})-\lambda_{l+p} \ \ {\rm si} \
1\leq l\leq r-p,$$
$$\beta'_l=\mu_{l+r-p+1} \ \ {\rm si} \ r-p+1\leq l\leq r+s-p-1.$$
On normalise ensuite, pour obtenir la suite
\m{(\alpha_1,\ldots,\alpha_{p+1},}\m{\beta_1,\ldots,\beta_{r+s-p-1})}
v\'erifiant
$$\sigg_{1\leq i\leq p+1}\alpha_i\dim(M^{(p)}_i) =
\sigg_{1\leq l\leq r+s-p-1}\beta_l\dim(N^{(p)}_l) = 1.$$
On a donc
$$\alpha_i = \q{\alpha'_i}{c}, \ \ \ \beta_l = \q{\beta'_l}{c},$$
avec
$$c=\sigg_{1\leq i\leq p}\lambda_im_i+\mu_1\biggl(
(\sigg_{p+1\leq j\leq r}m_j\dim(H_{1j}))-n_1\biggr).$$
On appelle
\m{(\alpha_1,\ldots,\alpha_{p+1},}\m{\beta_1,\ldots,\beta_{r+s-p-1})}
la {\em polarisation associ\'ee} \`a

\noindent
\m{(\lambda_1,\ldots,\lambda_r,}\m{\mu_1,\ldots,\mu_s)}. C'est une
polarisation de l'action de \m{G(p)} sur \m{W'(p)}. On supposera que les
\m{\alpha_i} et les \m{\beta_l} sont positifs.

\subsection{Comparaison des (semi-)stabilit\'es}

On veut comparer la (semi-)stabilit\'e d'un \'el\'ement de \m{W^0} avec
celle des \'el\'ements de \m{W'_0(p)} associ\'es.

\begin{xprop}
On suppose que
$$\sigg_{1\leq i\leq p}\lambda_im_i\ \leq \ \mu_1.$$
Si \ \m{w\in W^0} \ n'est pas \m{G}-(semi-)stable relativement \`a
\m{(\lambda_1,\ldots,\lambda_r,}\m{\mu_1,\ldots,\mu_s)}, alors \break
\m{z(w)\in W'_0(p)} \ n'est pas \m{G(p)}-(semi-)stable relativement \`a
\m{(\alpha_1,\ldots,\alpha_{p+1},}\m{\beta_1,\ldots,\beta_{r+s-p-1})}.
\end{xprop}

{\em\noindent D\'emonstration.} On ne traitera que le cas de la
semi-stabilit\'e, la stabilit\'e \'etant analogue. Posons \
\m{w=(\phi{li})_{1\leq i\leq r,1\leq l\leq s}}. Soient \
\m{M'_i\subset M_i}, \m{N'_l\subset N_l} \ des sous-espaces vectoriels
tels que
$$\epsilon=\sigg_{1\leq i\leq r}\lambda_i\dim(M'_i)-
\sigg_{1\leq l\leq s}\mu_l\dim(N'_l)>0$$
et \ \m{\phi_{li}(H^*_{li}\ot M'_i)\subset N'_l} \ pour \m{1\leq i\leq r},
\m{1\leq l\leq s}. On peut supposer que
$$N'_1=\sigg_{p+1\leq j\leq r}\phi_{1j}(H^*_{1j}\ot M'_j).$$
En effet, supposons que
$$k=\codim_{N'_1}(\sigg_{p+1\leq j\leq r}\phi_{1j}(H^*_{1j}\ot M'_j))>0.$$
On a, en posant \ \m{m'_i=\dim(M'_i)}, \m{n'_l=\dim(N'_l)},
\begin{eqnarray*}
\sigg_{p+1\leq i\leq r}\lambda_im'_i-\mu_1(n'_1-k) & = &
\epsilon-\sigg_{1\leq i\leq p}\lambda_im'_i+k\mu_1+
\sigg_{2\leq l\leq s}\mu_ln'_l\cr
 & > & k\mu_1 - \sigg_{1\leq i\leq p}\lambda_im'_i > 0\cr
\end{eqnarray*}
par hypoth\`ese. On peut donc au besoin remplacer
$$(M'_1,\ldots,M'_r,N'_1,\ldots,N'_s)$$
par
$$(0,\ldots,0,M'_{p+1},\ldots,M'_r,
\sigg_{p+1\leq j\leq r}\phi_{1j}(H^*_{1j}\ot M'_j),0,\ldots,0).$$
Soient \ \m{{M^{(p)}_i}'=M'_i\subset M^{(p)}_i} \ pour \m{1\leq i\leq p},
$${M^{(p)}_{p+1}}'=\ker(\sigg_{p+1\leq j\leq r}\phi_{1j})\cap
\biggl(\som_{p+1\leq j\leq r}(H_{1j}^*\ot M'_j)\biggr)\ \subset \
M^{(p)}_{p+1}=\ker(\sigg_{p+1\leq j\leq r}\phi_{1j}),$$
$${N^{(p)}_l}'=M'_{l+p}\subset N^{(p)}_l \ \ {\rm si} \ 1\leq l\leq r-p,
\ \ \ {N^{(p)}_l}'= N'_{l-r+p+1} \ \ {\rm si} \ r-p+1\leq l\leq r+s-p-1.$$
Il faut s'arranger pour trouver
$$z(w)=\pmatrix{\phi'_1 & \phi'_2 \cr x'_3 & x'_4}=
(\psi_{li})_{1\leq i\leq p+1,1\leq l\leq r+s-p-1}$$
de telle sorte que
$$\psi_{li}(H^{(p)*}_{li}\ot{M^{(p)}_i}')\subset {N^{(p)}_l}' \ \
{\rm pour} \ 1\leq i\leq p+1, 1\leq l\leq r+s-p-1.$$
On peut prendre \m{u} tel que
$$u\biggl(\som_{p+1\leq j\leq r,2\leq l\leq s}
(H_{1j}^*\ot B_{l1}^*\ot M'_j)\biggr)
\ \subset \ \som_{2\leq l\leq s}N'_l,$$
et \m{\alpha} tel que
$$\alpha\biggl(\som_{1\leq i\leq p}(H_{1i}^*\ot M'_i)\biggr) \ \subset
\som_{p+1\leq j\leq r}(H_{1j}^*\ot M'_j)$$
(car \ \m{N'_1 = \ov{\phi_2}(\som_{p+1\leq j\leq r}(H_{1j}^*\ot M'_j))} ).
Dans ce cas \m{z(w)} poss\`ede les propri\'et\'es voulues. \carre

\medskip

\begin{xcoro}
Si
$$\mu_1 \ \geq \ \q{1}{n_1+1}$$
et si \m{w\in W^0_p} est tel que \m{z(w)} n'est pas \m{G(p)}-(semi-)stable
relativement \`a

\noindent
\m{(\alpha_1,\ldots,\alpha_{p+1},}\m{\beta_1,\ldots,\beta_{r+s-p-1})},
alors \m{w} n'est pas \m{G}-(semi-)stable relativement \`a

\noindent\m{(\lambda_1,\ldots,\lambda_r,}\m{\mu_1,\ldots,\mu_s)}
\end{xcoro}

{\em\noindent D\'emonstration}. On applique la proposition pr\'ec\'edente
\`a la mutation inverse. \carre

\bigskip

Le r\'esultat suivant est imm\'ediat :

\begin{xprop}
Si
$$\mu_1 \ < \ \q{\sigg_{p+1\leq j\leq r}\lambda_jm_j}{n_1-1},$$
alors on a \ \m{W^{ss}\subset W^0_p}.
\end{xprop}

\medskip

\begin{xcoro}
Si
$$\mu_1 \ > \ \q{\sigg_{p+1\leq j\leq r}\lambda_jm_j}{n_1+1},$$
alors on a \ \m{W'(p)^{ss}\subset W'_0(p)}.
\end{xcoro}

\medskip

On en d\'eduit le

\begin{xtheo}
Si
$${\rm Max}(\q{1}{n_1+1},1-\sigg_{p+1\leq j\leq r}\lambda_jm_j) \
\leq \ \mu_1 < \q{\sigg_{p+1\leq j\leq r}\lambda_jm_j}{n_1-1},$$
alors il existe un bon quotient \ \m{W'(p)^{ss}//G(p)} \ si et seulement
si il existe un bon quotient \ \m{W^{ss}//G}, et dans ce cas les deux
quotients sont isomorphes.
\end{xtheo}

\bigskip

{\em\noindent Cas particuliers :}

\noindent 1 - Si \m{p=0}, la condition du th\'eor\`eme pr\'ec\'edent se
r\'eduit \`a
$$\mu_1 \ \geq \ \q{1}{n_1+1}.$$

\medskip

\noindent 2 - Si \m{s=1}, la condition du th\'eor\`eme pr\'ec\'edent se
r\'eduit \`a
$$\sigg_{p+1\leq j\leq r}\lambda_jm_j \ \geq \ \q{n_1-1}{n_1}.$$

\medskip

\noindent 3 - Si \m{p=0} et \m{s=1}, la condition du th\'eor\`eme
pr\'ec\'edent est toujours v\'erifi\'ee.

\bigskip

On suppose maintenant que \m{r=2} et \m{s=1}. Soient
$$\tau^* : H_{12}\ot A_{21}\lra H_{11}$$
la composition, et
$$\tau : H_{11}^*\ot A_{21}\lra H_{12}^*$$
l'application d\'eduite de \m{\tau}. On d\'eduit de ce qui pr\'ec\`ede
l'am\'elioration suivante du th\'eor\`eme 3.1 :

\begin{xtheo}
Si \ \m{r=2} \ et \ \m{s=1}, il existe un bon quotient projectif \
\m{W^{ss}//G} \ dans chacun des deux cas suivants :

\noindent 1 - On a
$$\q{\lambda_2}{\lambda_1}>\dim(A_{21}) \ \ \ {\rm et} \ \
\lambda_2\geq \q{\dim(A_{21})}{n_1} c(\tau,m_2).$$

\noindent 2 - On a
$$\lambda_1<\q{\dim(H_{11})}{n_1}, \ \ \lambda_2<\q{\dim(H_{12})}{n_1},
\ \ \lambda_2\dim(A_{21})-\lambda_1>
\q{\dim(A_{21})\dim(H_{12})-\dim(H_{11})}{n_1}, $$
et
$$\dim(H_{11})-\lambda_1n_1\geq c(\tau^*,m_1)\dim(A_{21}).$$
\end{xtheo}

\bigskip
\bigskip

\section{Applications}

\subsection{Vari\'et\'es de modules extr\'emales sur \proj{2}}

On note \m{Q} (resp. \m{Q_2}) le fibr\'e vectoriel sur \m{\proj{2}}
conoyau du morphisme canonique
$$\O(-1)\lra\O\ot H^0(\O(1))^* \ {\rm \ (resp. \ }
\O(-2)\lra\O\ot H^0(\O(2))^* \ {\rm \ )}.$$
Soient \m{m_1,m_2,n} des entiers positifs. On consid\`ere des morphismes
$$(*) \ \ \ \ (Q^*\ot\cx{m_1})\oplus(Q_2^*\ot\cx{m_2})\lra\O\ot\cx{n}.$$
Une polarisation est dans ce cas une suite \m{(\lambda_1,\lambda_2,\mu_1)}
de nombres rationnels positifs tels que
$$\lambda_1m_1+\lambda_2m_2 = \mu_1n = 1.$$
Cette polarisation est enti\`erement d\'etermin\'ee par le rapport
$$\rho = \q{\lambda_2}{\lambda_1}.$$
Supposons que \ \m{n < 3m_1+6m_2}. Alors les mutations des morphismes
pr\'ec\'edents sont du type
$$(**) \ \ \ \
\O\ot\cx{3m_1+6m_2-n}\lra(\O(1)\ot\cx{m_1})\oplus(\O(2)\ot\cx{m_2}).$$
La polarisation associ\'ee est \m{(\alpha_1,\beta_1,\beta_2)}, avec
$$\alpha_1=\q{1}{3m_1+6m_2-n}, \ \beta_1=\q{3-n\lambda_1}{3m_1+6m_2-n}, \
\beta_2=\q{6-n\lambda_2}{3m_1+6m_2-n}.$$
Elle est enti\`erement d\'etermin\'ee par le rapport
$$\rho' = \q{\beta_1}{\beta_2} = \q{3m_2\rho+3m_1-n}
{(6m_2-n)\rho+6m_1}.$$
Supposons que \ \m{m_1=1}. Alors on sait d'apr\`es \cite{dr_tr} construire
un bon quotient de l'ouvert des morphismes \m{G}-semi-stables de type
\m{(**)} d\`es que
$$\rho' > 3.$$

\subsubsection{Exemple 1}

Soit \m{m} un entier, \m{m\geq 0}. On consid\`ere des morphismes du
type
$$Q^*\oplus(Q_2^*\ot\cx{5m+2})\lra\O\ot\cx{29m+14}.$$
Si
$$\rho = \q{29}{12}+\epsilon,$$
avec \m{0<\epsilon\ll 1}, on sait d'apr\`es \cite{dr3} qu'un bon
quotient de l'ouvert des points \m{G}-semi-stables existe, et est
isomorphe \`a la vari\'et\'e de modules \m{M(4m+2, -2m-1, 2(m+1)^2)} des
faisceaux semi-stables au sens de Gieseker-Maruyama, de rang \m{4m+2} et de
classes de Chern \m{-2m-1} et \m{2(m+1)^2} sur \proj{2}. Dans \cite{dr3},
on ne donne pas de construction intrins\`eque du quotient
(on sait d\'ej\`a que \m{M(4m+2, -2m-1, 2(m+1)^2)} existe). Les
r\'esultats de \cite{dr_tr} ne permettent pas non plus de donner
directement l'existence du quotient. On peut cependant y parvenir en
utilisant le th\'eor\`eme 7.4, car dans ce cas
$$\rho' = 3 + \q{144\epsilon(m+1)}{29m+14+\epsilon(12m-24)} > 3.$$
On obtient des morphismes du type
$$\O\ot\cx{m+1}\lra\O(1)\oplus(\O(2)\ot\cx{5m+2}).$$

\subsubsection{Exemple 2}

Soit \m{m} un entier, \m{m\geq 0}. On consid\`ere des morphismes du
type
$$Q^*\oplus(Q_2^*\ot\cx{17m+8})\lra\O\ot\cx{99m+49}.$$
Si
$$\rho = \q{99}{41}+\epsilon,$$
avec \m{0<\epsilon\ll 1}, on sait d'apr\`es \cite{dr3} qu'un bon
quotient de l'ouvert des points \m{G}-semi-stables existe, et est
isomorphe \`a la vari\'et\'e de modules

\noindent
\m{M(7(2m+1), -4(2m+1), 32m^2+37m+11)}. Le th\'eor\`eme 7.4 permet de
construire le quotient, en utilisant les morphismes
$$\O\ot\cx{3m+2}\lra\O(1)\oplus(\O(2)\ot\cx{17m+8}).$$

\subsection{Un exemple sur \proj{n}}

On consid\`ere les morphismes
$$(\phi_1,\phi_2) : \O(-2)\oplus\O(-1)\lra\O\ot\cx{n+2}$$
sur \proj{n}. Une polarisation est dans ce cas un triplet
\m{(\lambda_1,\lambda_2,\mu_1)} de nombres rationnels positifs tel que
$$\mu_1=\q{1}{n+2}, \ \lambda_1+\lambda_2 = 1.$$
Il revient au m\^eme de se donner
$$\rho = \q{\lambda_2}{\lambda_1}.$$
On sait construire des bons quotients (en utilisant les r\'esultats de
\cite{dr_tr})) d\`es que
$$\rho > n+1.$$
Mais dans ce cas le quotient est vide ! En effet, il existe toujours un
sous-espace vectoriel \ \m{H\subset\cx{n+2}} \ de dimension \m{n+1} tel
que \ \m{\imm(\phi_2)\subset \O\ot H}. On doit donc avoir, si
\m{(\phi_1,\phi_2)}
est \m{G}-semi-stable relativement \`a \m{(\lambda_1,\lambda_2,\mu_1)},
$$\lambda_2 - \q{n+1}{n+2} \leq 0,$$
c'est-\`a-dire \ \m{\rho\leq n+1}.

On emploie maintenant le th\'eor\`eme 7.4, et on consid\`ere donc les
morphismes
$$\O\ot\cx{\q{n(n+3)}{2}}\lra Q_2\oplus Q,$$
o\`u, comme dans le \para 8.1, \m{Q} (resp. \m{Q_2}) d\'esigne le
le fibr\'e vectoriel sur \m{\proj{n}} conoyau du morphisme canonique
$$\O(-1)\lra\O\ot H^0(\O(1))^* \ {\rm \ (resp. \ }
\O(-2)\lra\O\ot H^0(\O(2))^* \ {\rm \ )}.$$
En utilisant les r\'esultats de \cite{dr_tr} on parvient \`a construire
un bon quotient projectif d\`es que
$$\rho \geq 1 - \q{2n}{(n+1)(n+4)}.$$
Les valeurs {\em singuli\`eres} de \m{\rho} sont par d\'efinition celles
pour lesquelles la \m{G}-semi-stabilit\'e n'implique pas la
\m{G}-stabilit\'e. Ces valeurs sont exactement les nombres
$$\rho_k = \q{k}{n+2-k}$$
pour \m{1\leq k\leq n+1} . Dans ce cas un morphisme \m{(\phi_1,\phi_2)}
$G$-semi-stable non $G$-stable est construit de la fa\c con suivante :
on consid\`ere un sous-espace vectoriel \ \m{H\subset\cx{n+2}} \ de
dimension \m{k}, et on prend pour \m{\phi_2} un morphisme tel que
\ \m{\imm(\phi_2)\subset\O\ot H} \ et que $H$ soit le plus petit
sous-espace vectoriel ayant cette propri\'et\'e. On prend pour \m{\phi_1}
un morphisme tel que l'application lin\'eaire induite
$$H^0(\O(2))^*\lra\cx{n+2}$$
soit surjective.

On obtient donc au total \ \m{2\lbrack\q{n}{2}\rbrack +2} \ quotients
distincts et non vides, dont \ \m{\lbrack\q{n}{2}\rbrack} \ sont
singuliers. Ils sont de dimension \
\m{\q{(n+2)(n^2+3n-2)}{2}}, sauf celui correspondant \`a \m{\rho_{n+1}},
qui est de dimension \m{\q{n(n+3)}{2}}.

\bigskip

On peut g\'en\'eraliser ce qui pr\'ec\`ede et obtenir des bons quotients
projectifs d'espaces de morphismes du type
$$\O(-p-q)\oplus\O(-p)\lra\O\ot\cx{n+2}$$
sur \proj{n}, \m{p,q} \'etant des entiers positifs.


\bigskip
\bigskip

\begin{thebibliography}{99}
\bibitem{dr1}Dr\'ezet, J.-M. {\em Fibr\'es exceptionnels et suite spectrale
de Beilinson g\'en\'eralis\'ee sur $\proj{2}(\cx{})$} .
Math. Ann. 275 (1986), 25-48.
\bibitem{dr2}Dr\'ezet, J.-M. {\em Fibr\'es exceptionnels et vari\'et\'es de
modules de faisceaux semi-stables sur \proj{2}($\hskip -1pt\cx{} \hskip
1pt)$.} Journ. Reine Angew. Math. 380 (1987), 14-58.
\bibitem{dr3}Dr\'ezet, J.-M. {\em Vari\'et\'es de modules extr\'emales de
faisceaux semi-stables sur $\proj{2}(\cx{})$ .}Math. Ann. 290 (1991),
727-770.
\bibitem{dr_lp}Dr\'ezet, J.-M., Le Potier, J. {\em Fibr\'es stables et fibr\'es
exceptionnels sur \proj{2}}. Ann. Ec. Norm. Sup. 18 (1985), 193-244.
\bibitem{dr_tr}Dr\'ezet, J.-M., Trautmann, G. {\em Moduli spaces of morphisms
of sheaves and quotients by non-reductive groups} Preprint (1995).
\bibitem{go_ru} Gorodentsev, A.L., Rudakov, A.N. {\em Exceptional vector
bundles on projective spaces .} Duke Math. Journ. 54 (1987), 115-130.
\bibitem{karpov}Karpov B.V. {\em Semi-stable sheaves on a two-dimensional
quadric and Kronecker modules.} Math. Izvestiya AMS transl. 40 (1993),
33-66.
\bibitem{king} King, A. {\em Moduli of representations of finite dimensional
algebras .} To appear in the Quarterly Journal of Mathematics.
\bibitem{miro} Mir\'o-Roig, R.M. {\em Some moduli spaces for rank 2 stable
reflexive sheaves on \m{\proj{3}}. } Trans. Amer. Math. Soc. 299 (1987),
699-717
\bibitem{miro-trm}Mir\'o-Roig, R.M., Trautmann, G.{\em The moduli scheme
\m{M(0,2,4)} over \m{\proj{3}}}. Math. Z. 216 (1994), 283-315.
\bibitem{mumf}Mumford, D., Fogarty, J. {\em Geometric invariant theory.} Ergeb.
Math. Grenzgeb. Bd. 34. Berlin Heidelberg New-York : Springer (1982)
\bibitem{news}Newstead, P.E. {\em Introduction to moduli problems and orbit
spaces. } TIFR Lect. Notes. Math. vol. 51. Berlin Heidelberg New-York :
Springer (1978)
\bibitem{oko}Okonek, C. {\em Moduli extremer reflexiver Garben auf
\m{\proj{n}}. } Journ. Reine Angew. Math. 338 (1983), 183-194
\end{thebibliography}
\end{document}